\documentclass{amsart}
\usepackage{lmodern}
\usepackage{hyperref}
\usepackage[scale=2]{ccicons}
\usepackage[frenchb,english]{babel}
\usepackage[T1]{fontenc}
\usepackage{multirow}
\usepackage{adjustbox}
\usepackage{amsmath}
\usepackage{bbm}
\numberwithin{equation}{section}
\usepackage{amsfonts}
\usepackage{amssymb}
\usepackage{hyperref}
\usepackage{dsfont}
\usepackage{caption}
\usepackage{tikz}
\usepackage{amsthm}
\usepackage{enumitem} 
\usepackage{color}
\usepackage{graphicx}
\usepackage{hyperref}
\usepackage{tikz}
\usetikzlibrary{automata,topaths}
\usetikzlibrary{arrows,positioning,shadows,fit,shapes}
\usepackage{enumitem}
\usepackage{mathtools}

\newtheorem*{remark}{Remark}

\newtheorem{theorem}{Theorem}[section]
\newtheorem{lemma}{Lemma}[section]
\newtheorem{proposition}{Proposition}[section]
\newtheorem{corollary}{Corollary}[section]
\newtheorem{definition}{Definition}[section]
\newcommand{\R}{\mathbb{R}}
\newcommand{\N}{\mathbb{N}}

\newcommand{\Z}{\mathbb{Z}}

\renewcommand{\P}{\mathbb{P}}
\newcommand{\Q}{\mathbb{Q}}

\begin{document}

\title[Stable web]{Construction and convergence results for stable webs}

\author{Thomas Mountford}

\address{Thomas Mountford\\\'Ecole Polytechnique F\'ed\'eral de Lausanne\\Department of Mathematics\\EPFL SB MATH PRST\\
MA B1 517 (Bâtiment MA)
Station 8
CH-1015 Lausanne\\
Switzerland}
\email{thomas.mountford@epfl.ch}
\urladdr{http://people.epfl.ch/thomas.mountford}

\author{Krishnamurthi Ravishankar}

\address{Krishnamurthi Ravishankar\\NYU-ECNU, Institute of Mathematical Sciences at NYU-Shanghai,
  Shanghai 200062, China }
\email{ravi101048@gmail.com}
\urladdr{}
\thanks{This work was partially supported by  Simons Foundation Collaboration grant 281207 awarded to K. \,Ravishankar. We thank NYU-Shanghai, where initial part of the work was done for their hospitality and support.}


\subjclass[2010]{Primary 60K35; Secondary 60F17, 60J15}

\keywords{Coalescing random walk, Stable processes, Brownian motion }

\maketitle

\begin{abstract}
    We introduce a new metric for collections of aged paths and a robust set of criteria for compactness for a set of collection of aged paths in the topology corresponding to this metric. We show that the distribution of stable webs ($1< \alpha \leq 2$)  made up of collections of stable paths is tight in this topology.  We then show the weak convergence of appropriately normalized systems of coalescing random walks in the domain of attraction of stable laws for $1 < \alpha \leq 2$ under this metric to the corresponding stable web. We obtain some path results in the brownian case.

\end{abstract}


\section{Introduction}
In this article we consider collections of aged paths following \cite{MRV}.  Our motivations in that article 
were (and remain) to produce a topology for a system of coalescing stable processes (of index $\alpha > 1$) so that,  equipped with an associated age process, systems of coalescing random walks in the domain of attraction of the given stable process, when suitably renormalized, converge  to the collection of stable aged paths.

Our chief motivation for introducing aged paths was that the (now) classical approach of \cite{FINR} for the convergence to the Brownian web did not seem to permit analogous results for convergence to systems of 
coalescing stable processes or the construction of such a system.  A big obstacle was the fact that
for, say, the space time set $[0,1]^2$ a  system of coalescing stable processes ($1 <\alpha <2$)  beginning at a countable dense set of
space time points did not have a maximum displacement at times $t > 0$.  (We do not consider $\alpha \leq 1$ since, as is well known, in this case two independent processes starting from distinct points never meet a.s..) This behaviour for $1< \alpha < 2$ is  in contrast to that for Brownian paths.   Nonetheless our system is well defined for Brownian paths resulting in a distinct topology to that for the web (even if we put aside the addition of the age functions).

In this work, we change slightly the definition of aged paths and introduce a new metric.  
This metric is very similar to that of \cite{MRV} but we believe it superior  since  it is homogenous in space and time rather than privileging in importance path behaviour at points where the time or spatial value is an integer.  
We believe, in that it removes some conditions for a set of collections  of paths to be compact, that it also
constitutes a simplification.  We emphasize
that though our motivation is the study of coalescing systems in this part,  Section 2, coalescence need never be mentioned.
This section concludes with a robust set of criteria for compactness of sets of
aged path collections, Theorem \ref{thm3}.

These compactness criteria are then applied in Section 3 to establish Theorem \ref{thm4}, the result that
collections of  coalescing 
random walks in the domain of attraction of stable laws converge in distribution to the corresponding stable web.
This follows a similar result in \cite{MRV} but we emphasize that the result is new in that it is established for a different topology and this leads to a simpler proof.

The treatment also includes $\alpha = 2$, that is Brownian motion which was outside the criteria furnished in the earlier paper.

 A good property of  the topology proposed here is that it permits 
convergence of (suitably renormalized) systems of coalescing random walks to the "Brownian" limit if the walks have kernel with mean zero and finite variance.  That is the conditions of Donsker's Invariance principle.
This is in sharp contrast to the situation for the Brownian web (see \cite{NRS} and \cite{BMRV}) where for convergence the kernel must have  "loosely" 3'rd moment. For the case $1 < \alpha < 2$ we do not attempt to give a proof with the widest possible class of random walks, contenting our selves with symmetric cases which are very regular.

In the last section we work exclusively with the the stable web for $\alpha = 2$ (which we call the {\it Brownian stable web}).  In this case our paths are continuous and so  spatial order is preserved
 and we can make the argument that (with high probabiliy) well separated space time regions are ``ìnsulated"
from one another and have almost independence properties.  This and the very thin tails for Gaussian
distributions enable us to prove (Proposition \ref{prop: brown 1}) that for all our aged paths we can extend the domain of defintion from
$(\sigma, \infty) $ to $[\sigma, \infty) $ in a continuous fashion.  Here $ \sigma$  is the birth time of an aged path. We would like to do this for more general stable webs but are unable to do so at the moment.
 We also prove in this section that distinct paths cannot have the same $\sigma$ (Proposition \ref{prop: brown 2}).  The underlying intuition suggests again that this brownian result   holds more widely.
\section{The construction}
In this section, we introduce the space and metric that we will be working with and show some simple properties. We first wish to define the space of \textit{locally finite} sets of aged paths. Our definitions will be similar to those given in \cite{MRV}  but  with slight differences.
\begin{definition}
\label{def: age path}
An aged path will be a pair of c\`ad\`ag functions $(\gamma,a)$ defined on an open interval $(\sigma,\infty)$, so that in addition
\begin{enumerate}[label=(\roman*)]
\item $\lim_{s\downarrow \sigma}a(s)=0$
\item $a(t+s)-a(s)\geq t,\forall t>0,s>\sigma$
\item $\gamma,a$ do not jump at the same time.
\end{enumerate}
\end{definition}
\begin{remark}
In \cite{MRV}, the age does not satisfy (i), but only $\lim_{s\downarrow\sigma}a(s)\geq0$;  equally in \cite{MRV} there is no condition analagous to (iii).
\end{remark}
For an aged path $(\gamma,a)$ and $t\geq 1$, the operation $\Pi_t$ acts on $(\gamma,a)$ by $\Pi_t(\gamma,a)=(\gamma_t,a_t)$ where 
 $\gamma_t$ and $a_t$ are simply restrictions of $\gamma,a$ to $[b_t,t]$, where 
$$b_t=b_t(\gamma,a):=\inf\{-t\leq s\leq t:\gamma(s)\in [-t,t],a(s)\geq 2^{-t}\}.$$
Notice that the result of this  projection, $(\gamma_t,a_t)$ is not an aged path since it violates (i). Given a collection of aged paths $\Gamma=\{(\gamma_i,a_i),i\in I\}$, $\Pi_t(\Gamma)$ is equal to

$$\left\{\Pi_t(\gamma,a): (\gamma,a) \in \Gamma, \nexists (\gamma',a')\in\Gamma\text{ so that }(\gamma',a')=(\gamma,a)\text{ on }[b_t(\gamma,a),t]\text{ but }b_t(\gamma',a')<b_t(\gamma,a)\right\}\backslash \varnothing,$$

where $\varnothing$ represents the trivial aged path with empty domain. $\Pi_t(\Gamma)$ is simply the collection of \textit{maximal} elements of $\Pi_t(\gamma,a)$ with $(\gamma,a)\in\Gamma$. 
 Note the definition presupposes the existence of such maximal paths.     
\begin{definition}
A collection of aged paths $\Gamma$ is said to be locally finite if $\forall t\geq 1$, $|\Pi_t(\Gamma)|<\infty$.
\end{definition}
Let $\mathcal{X}$ be the space of collections of locally finite  aged paths.   $\mathcal{X}^\prime \supset \  \mathcal{X}$
will denote the space of locally finite paths saitsfying (ii) and (iii) but not necesarily (i),  though it is still required that $\lim_{s \downarrow \sigma} a(s) \geq \ 0$.

As in \cite{MRV}, we use the following small generalization of the Skorohod metric. For c\`adl\`ag functions
$$f:[a,b]\rightarrow \R\text{ and }g:[c,d]\rightarrow\R.$$
The distance between $f$ and $g$ is
\begin{equation}
\label{function metric}
d(f,g)=\inf_{\tau:[a,b]\rightarrow[c,d]}\sup_{s\in[a,b]}|\tau(s)-s|+|f(s)-g(\tau(s))|,
\end{equation}
where the infimum is over all increasing homeomorphism $\tau:[a,b]\rightarrow[c,d]$.

We define the metric between two locally finite collections of aged paths as
$$d_{\mathcal{X}}(\Gamma,\Gamma')=\int_1^{\infty}e^{-t}\left(d(\Pi_t(\Gamma),\Pi_t(\Gamma'))\wedge 1\right) dt,$$
where (abusing notation and extending $d$ to be a metric between finite sets of paths) $d(\Pi_t(\Gamma),\Pi_t(\Gamma'))$ is the Hausdorff metric  $d(\Pi_t(\Gamma), \Pi_t(\Gamma^{'}))=$

\begin{align*}
\max \bigg\{ \sup_{(\gamma,a)\in\Pi_t(\Gamma)}&\inf_{(\gamma ^{'},a^{'})\in\Pi_t(\Gamma')}d((\gamma,a),(\gamma^{'},a^{'})),   \\   
& \sup_{(\gamma ^{'},a^{'})\in\Pi_t(\Gamma^{'})} \inf_{(\gamma,a)\in\Pi_t(\Gamma)} d((\gamma,a),(\gamma ^{'},a^{'}))\bigg\}  , \\
\end{align*}

where $d((\gamma,a),(\gamma ^{'},a^{'}))=d(\gamma,\gamma^{'})\vee d(a,a^{'})$. 

We will use $d_{\mathcal{X}}$  to also denote this metric on the space $\mathcal{X}^\prime $.

Notice that instead of using $e^{-t}$ in the metric $d_{\mathcal{X}}$, it can be replaced by any decreasing, positive and integrable function $g$ on $[0,\infty)$. It is easy to see that $\forall \epsilon>0$, $\exists \delta>0$ such that
$$\int_1^{\infty}e^{-t}\left(d(\Pi_t(\Gamma),\Pi_t(\Gamma'))\wedge 1\right)dt<\delta$$
implies
$$\int_1^{\infty}g(t)\left(d(\Pi_t(\Gamma),\Pi_t(\Gamma'))\wedge 1\right)dt<\epsilon,$$
and vice versa.   

So the choice of function $e^{-t}$ does not change the topology. Another choice in defining the metric is $2^{-t}$ in the restriction $\Pi_t$. We could have replaced $2^{-t}$ by any continuous decreasing function $h:[1,\infty)\rightarrow (0,\infty)$ that tends to $0$ at $\infty$. We could also define $\Pi_t^h(\gamma,a)$ as the restriction of $\gamma$ and $a$ to $[b^h_t,t]$, where
$$b^h_t=\inf\{s\geq -t:\gamma(s)\in[-t,t],a(s)\geq h(t)\}.$$
We could then have defined
$$d^h_{\mathcal{X}}(\Gamma,\Gamma')=\int_1^{\infty}e^{-t}\left(d(\Pi^h_t(\Gamma),\Pi^h_t(\Gamma'))\wedge 1\right)dt$$
as before. To show that $d^h_{\mathcal{X}}$ gives the same topology requires more work. We need the following lemma whose content is well known  and whose proof is left to the reader.

\begin{lemma}
\label{lem: cadlag}
For c\`adl\`ag functions $f:[a,b]\rightarrow\R$ and $f_n:[a_n,b_n]\rightarrow\R$ satisfying $d(f,f_n)\rightarrow 0$ as $n\rightarrow\infty$, we have for any $a\leq t\leq b$,
\begin{enumerate}[label=(\roman*)]
\item if $t$ is a continuity point of $f$, then $\forall\{t_n\}_{n=1}^{\infty}$ such that $t_n\in[a_n,b_n]$ and $t_n\rightarrow t$, we have
$$f_n\mid_{[t_n,b_n]}\overset{d}{\rightarrow}f\mid_{[t,b]}.$$
\item if $t$ is a jump point of $f$ (so $t\neq a$), then $\exists \{t_n\}_{n=1}^{\infty}$ such that $t_n\in[a_n,b_n]$ and $t_n\rightarrow t$, we have
$$f_n(t_n)\rightarrow f(t),f_n(t_n-)\rightarrow f(t-)\text{ and }f_n\mid_{[t_n,b_n]}\overset{d}{\rightarrow}f\mid_{[t,b]}.$$
\end{enumerate}
\end{lemma}

\begin{remark}
It is not true in general that $\forall t\in [a,b]$ and $t_n\rightarrow t$, we have
$$f_n\mid_{[t_n,b_n]}\overset{d}{\rightarrow}f\mid_{[t,b]}.$$
For example when $t$ is a discontinuity point in $(a,b]$ and $\forall n  \ f_n \equiv f $ and $t_n<t$.
\end{remark}

\begin{remark}
 The sequence $ \{t_n\} $ in part (ii) is essentially unique in that  if $\{t_n^\prime \}$ is another such sequence, then
$ t_n = t_n ^\prime $ for $n$ large.
\end{remark}
\begin{corollary} \label{nice}
For $f_n, \ f$ as above, if $u,v$ are continuity points of $f$ in $(a,b)$, then $f_n \mid_{[u,v]}
\overset{d}{\rightarrow}f\mid_{[u,v]}$ 
\end{corollary}
\begin{remark}
In this corollary, $f_n$ need not be defined on $[u,v]$ for every $n$ but will be so defined for $n$ sufficiently large.
\end{remark}

\begin{proposition}
\label{prop:same top}
$d^h_{\mathcal{X}}$ gives the same topology as $d_{\mathcal{X}}$ when $h:[1,\infty)\rightarrow (0,\infty)$ is continuous, decreasing and satisfies $\lim_{t\rightarrow\infty}h(t)=0$.
\end{proposition}

\begin{proof}

We wish to show that convergence in one metric gives convergence in the other.  We show that convergence in $d^h_{\mathcal{X}}$ implies $d_{\mathcal{X}}$ convergence.   The converse implication follows in the same manner To do this it is only necessary to 
consider the behaviour for $t$ in a compact interval.  

We first fix interval $[N,N+1)$ and suppose $\Gamma_n\overset{d^h_{\mathcal{X}}}{\rightarrow}\Gamma$. We can choose $m>N$ large enough so that
$$h(m)<2^{-(N+1)}.$$
By taking subsequences, we can assume that
$$d(\Pi^h_t(\Gamma_n),\Pi^h_t(\Gamma))\rightarrow 0 \text{ for a.e. }t\in [m,m+1)\text{ as }n\rightarrow\infty.$$
Fix $T\in [m,m+1)$ so that we have $d(\Pi^h_T(\Gamma_n),\Pi^h_T(\Gamma))\rightarrow 0$ as $n\rightarrow\infty$. We wish to show that
$$d(\Pi_t(\Gamma_n),\Pi_t(\Gamma))\rightarrow 0\text{ for a.e. }t\in [N,N+1)\text{ as }n\rightarrow\infty.$$
As $N$ can be  arbitrarily large this will suffice to show $d_{\mathcal{X}}$ convergence.

Take $V_N$ to be  the  subset of $[N,N+1)$ so that one of the following conditions holds: for $t \in V_N$:
\begin{enumerate}[label=(\roman*)]
\item $-t$, $t$ or $2^{-t}$ is a discontinuity point for $(\gamma,a)\in\Pi_T(\Gamma)$ or for $t\mapsto b_t(\gamma,a)$
\item $-t$ or $t$ is a \textit{jump value} for $\gamma$ or $a$ with $(\gamma,a)\in\Pi_T(\Gamma)$. A c\`adl\`ag  function $f$ has a \textit{jump value} of $c$ if $\exists s\in[-,t,t]$ such that $f(s)\neq f(s-)$ and $f(s)=c$ or $f(s-)=c$.
Further  $\gamma(s) \ne t$ for $s$ a jump time of $a$ and vice versa. 
\item $\gamma$ is discontinuous at $s=\inf\{u:a(u)\geq 2^{-t}\}$
\item $t$ is a \text{``1corner" point} of $(\gamma,a)\in\Pi_T(\Gamma)$:
\begin{align*}
b_t=t&\text{ and either }\gamma(t)=t\text{ or }-t\\
b_t=t&\text{ and }a(t)=2^{-t}
\end{align*}
\item $s\mapsto|\Pi_s(\Gamma)|$ jumps at $t$.
\end{enumerate}
It
is easily seen that $V_N$ is countable and so of Lebesgue measure zero.
 By condition (v) defining $V_N$ it is only necessary to show that for  $(\gamma,a)\in \Pi_t(\Gamma)$, we have appropriate convergence for functions in  $(\gamma,a)\in \Pi_t(\Gamma_n )$.
For $t\notin V_N, N\leq t<N+1$, we have that for $(\gamma,a)\in \Pi_t(\Gamma)$, there is $(\gamma',a')\in \Pi_T(\Gamma)$ so that $\Pi_t(\gamma',a')=(\gamma,a)$. By hypothesis, there exists $(\gamma'_n,a'_n)\in\Pi_T(\Gamma_n)$ with $d((\gamma'_n,a'_n),(\gamma',a'))\rightarrow 0$ as $n\rightarrow\infty$. Let $(\gamma_n,a_n)=\Pi_t(\gamma'_n,a'_n), b_t=b_t(\gamma,a)$ and $b_{n,t}=b_t(\gamma_n,a_n)$. We want to show that $d(\gamma_n,\gamma)$ and $d(a_n,a)$ tend to zero as $n$ tends to infinity. We only treat $\gamma$ and the proof for $a$ is similar. 

There are only four (not mutually exclusive) possibilities :
\begin{enumerate}
\item $\exists \epsilon>0$ so that $\gamma_s\notin [-t,t]$ for $s\in(b_t-\epsilon,b_t)$ and $\gamma$ is continuous at $b_t$.
\item  $\exists \epsilon>0$ so that $\gamma_s\notin [-t,t]$ for $s\in(b_t-\epsilon,b_t)$ and $\gamma$ is discontinuous at $b_t$.
\item $a(s)<2^{-t},\forall s<b_t$ and $a$ is continuous at $b_t$.
\item $a(s)<2^{-t},\forall s<b_t$ and $a$ is discontinuous at $b_t$.
\end{enumerate}
We first consider (2). Hypothesis (iii) in Definition \ref{def: age path} ensures that $a$ is continuous at $b_t$ and condition (iii) of $V_N$ ensures that $a(b_t)>2^{-t}$. At $b_t$, $\gamma$ jumps from a point in $[-t,t]^c$ to a point in $(-t,t)$  by (ii) above.   So we can apply Lemma \ref{lem: cadlag} to see that there exists $c_{n,t}  \rightarrow b_t $
so that
there is a corresponding jump of $\gamma_n$ at $c_{n,t}$ and that  
$$\gamma_n\mid_{[c_{n,t},t]}\overset{d}{\rightarrow}\gamma\mid_{[b_t,t]}.$$
It remains to check that $c_{n,t}=b_{n,t}$ for $n$ large enough.   It is immediate that $c_{n,t} \geq b_{n,t}$ for $n$ large.   We have that for some $\epsilon^\prime >0$, $$\gamma_n(s)\notin [-t,t],\text{ for }s\in(c_{n,t}-\epsilon^\prime,c_{n,t}).$$
So if $c_{n,t}\neq b_{n,t}$ for $n$ arbitrarily large, we must have for some subsequence of $\{b_{n,t}\}_{n=1}^{\infty}$,
$$\lim_{k\rightarrow\infty}b_{n_k,t}=b'\leq b_t-\epsilon^\prime .$$
There are two possibilities: $a(b')<2^{-t}$ or $\gamma(b')\notin [-t,t]$. The first possibility contradicts with the fact that $a_{n_k}(b_{n_k,t})\geq 2^{-t}$. The second possibility requires that $t$ is a \textit{jump value} for $\gamma$ which contradicts with condition (ii) of $V_N$ because of the following lemma  whose proof is again left to the reader.
\begin{lemma}
Suppose $t_0$ is a discontinuity point of a c\`adl\`ag function $f$ and $f_n\overset{d}{\rightarrow}f,t_n\rightarrow\ t_0$. If $\lim_{n\rightarrow\infty}f_n(t_n)=c$, then $$f(t_0)=c\text{ or }f(t_0-)=c,$$which means that $c$ is a \textit{jump value} of $f$.
\end{lemma}
Next, we consider  (4). Again if $a$ is discontinuous at $b_t$, we must have that at $b_t$, $a$ jumps from a point in $(0,2^{-t})$ to a point in $(2^{-t},\infty)$   and $\gamma (b_t) \in (-t,t)$ by (ii) above. Thus we have $b_{n,t}\rightarrow b_t$ and $\gamma_n\overset{d}{\rightarrow}\gamma$ by Lemma \ref{lem: cadlag}. For case  (1), we argue again as in (2) that $b_{n,t}\rightarrow b_t$ and Lemma \ref{lem: cadlag} suffices  as $b$ is continuous at $t$. For case  (3),  using that $b$ is continuous at $t$, $b_{n,t}\rightarrow b_t$ so the only problem arises if $\gamma$ has a jump at $b_t$ but condition (iii) of $V_N$ removes this possibility.
\end{proof}

For  increasing right continuous $M_t \in [1,\infty) $ and  $\delta_t(n)$ increasing in $t$ but tending to zero as $n$ tends to infinity for each $t$, let $\mathcal{K}$ (= $\mathcal{K}(M, \delta) $ ) be the set of collections of aged paths in $\mathcal{X}$ satisfying the following conditions:
\begin{enumerate}[label=(\Alph*)]
\item   $|\Pi_tK|\leq M_t$ for any $K\in\mathcal{K}$.\label{i} 
\item $\forall (b_t,\gamma,a)\in\Pi_tK$, $a(s)\leq M_t$ and $|\gamma(s)|\leq M_t$ for $  b_t\  \leq s\leq t$ and any $K\in\mathcal{K}$.\label{ii}
\item (Recall for a cadlag function $f$ defined on interval $I \ = \ [a,b], \ \omega (\delta,f,I) \ = \ \inf_{n, t^n_i} \sup_i\{|f(s)-f(t)|: s,t, \in [t^n_i,t^n_{n+1})\}$
where the infimum is taken over 
positive integers $n$ and $a=t^n_0 < t^n_1 \cdots t^n_n =b$ with  $t^n_i - t^n_{i-1} \geq \delta \ \forall i $ (with the last interval taken to be closed ).  See e.g. \cite{ek}).
$\forall n \geq 1 $ $\forall (b_t,\gamma,a)\in\Pi_t(K)$ has $\omega(2^{-n},\gamma^t,[-(t+1),t+1])\vee\omega(2^{-n},a^t,[-(t+1),t+1])\leq \delta_t(n)$ for any $K\in\mathcal{K}$, where $\gamma^t $ is the function 
defined by 

$\gamma^t(s)  = \gamma(b_t)$ for  $s \leq b_t ; $= $\gamma(s)$ for $t \geq s \geq b_t; \ = \ \gamma (t) $ for 
$t \leq s \leq t+1$ and similarly for $a^t$, though respecting the at least linear growth stipulation.\label{iii}
\item $\forall (b_t,\gamma,a)\in\Pi_tK$, and $s',s\in[b_t,t]$ such that $\gamma$ jumps at $s$ and $a$ jumps at $s'$ with sizes $|\gamma(s)-\gamma(s-)|\geq 2^{-n}$ and $a(s')-a(s'-)\geq 2^{-n}$, then $|s-s'|\geq \delta_t(n)$ for any $K\in\mathcal{K}$.\label{iv}

\item 
If $\Pi_t(\gamma,a)\neq \varnothing$ for $(\gamma,a) \in K$, then  
 there exist $(\gamma',a') $ and $s > \sigma '$ such that 
$\Pi_t(\gamma',a') \ = \ \Pi_t(\gamma,a)$,  $|\gamma'_s|\leq M_t$  and  so that $a'(s) \in ( 2^{-2t},2^{-3t/2})$\label{v}

\end{enumerate}

 Remark: We note that conditions \ref{ii} and \ref{iii}  correspond to  standard conditions for compactness of paths of c\`adl\`ag functions. Condition \ref{i} is needed to ensure that sub sequential limits remain in $\mathcal{X}$. Conditions \ref{iv}  is needed to ensure compactness of the aged paths and \ref{v} is used to extend the convergence from finite space-time region with positive lower bound on age to the whole space.
\begin{theorem} \label{thm3}
The set $(\mathcal{K},d_{\mathcal{X}})$  defined above is compact.
\end{theorem}

\begin{proof}

To show $\mathcal{K}$ is compact, it is sufficient to show that for any sequence of collections $\{K_n\}_{n\in\N},  K_n \in \mathcal{K}$, we can find a subsequence that converges to an element $K\in\mathcal{K}$. For any integer $N > 1$ and $K \in \mathcal{K}$,
 $\text{Range}(\Pi_{N}(\gamma,a))\subset [-M_{N},M_{N}]\times[2^{-N},M_{N}]$ for any $(\gamma,a)\in K$.  By property \ref{i}, we know that there  are  at most $M_N$ paths in $\Pi_{N}K$ for any $K\in\mathcal{K}$. For any element $(\gamma_N,a_N)\in\Pi_{N}K$, its range is a subset of a compact set in $\R^2$ as observed above.  In order to avoid 
starting point anomalies  (since we are dealing with $\Pi_t$ for all $t \in [1,\infty)$), we regard an element $(\gamma_N, a_N) $ in $\Pi_{N}K$
as actually having domain
$[-(N+1),N+1]$,  by extending the functions as follows; $\gamma^t$ is taken to be constant and $a^t$ linear with slope one on $[-(N+1),b_N]$ and on $[N,N+1]$ following \ref{iii} above .
Theorem 3.6.3 of \cite{ek} states an equivalent condition when $A\subset D_E[0,\infty)$ is compact, where $D_E[0,\infty)$ is the set of c\`adl\`ag functions valued in a complete metric space $(E,r)$. This theorem is also true when the time set is a finite interval $[a,b]$ instead of $[0,\infty)$.
Hence thanks to property \ref{ii} and \ref{iii} and an adaption of Theorem 3.6.3 of \cite{ek}, $\{ ( \gamma_N,a_N) \in \Pi_{N}K $ for some $K \in \mathcal{K}\} := \mathcal{A}$ is in a compact subset of  $D([-(N+1), N+1])$. for any $N$, where the  $D(I)$ are the sets of c\`adl\`ag functions defined 
on  compact intervals $I$ equipped with metric (\ref{function metric}) which is equivalent to the Skorohod metric on compact sub-interval of $\R$.  
Since using Theorem 3.6.3 of \cite{ek} it follows that $\Pi_N K$ is a compact subset of $\mathcal{A}$ for any $K \in \mathcal{K}$, we obtain that $\cup_{K \in \mathcal{K}} \Pi_{N}K$ is in a compact subset of $(\mathcal{X}^\prime ,d_{\mathcal{X}})$. 
We can find a subsequence $\{K_{n_j^N}\}_{j\in\N}$ of $\{K_n\}_{n\in\N}$ such that
 $\Pi_{N}K_{n_j^N}$ converges for $ N \geq 1$, 
to some collection $ \tilde{K}(N) \in \mathcal{X}^\prime$  with $| \tilde{K}(N)| \leq M_N$. 
By the usual diagonal procedure we obtain a subsequence $K_{m_j}$ so that for each $N$
$$
\Pi_{N}K_{m_j} \  \rightarrow \ K'(N)
$$
In the following, to avoid notational encumbrance, we regard the subsequence as being the original sequence,
that is  we write $K_j$ for $ K_{m_j} $
We now define for $t \geq 1, K(t) = \Pi_t K'(N)$ for $N > t$.  Corollary \ref{nice} implies that this definition is consistent.
We note that with this definition it need not be the case that $K(N) \ = \ K'(N)$.  Consider for simplicity
aged path collections consisting of a single function: $K_n \ = \ \{(\gamma_n, a_n)\}$ where
$$
a_n(s) = s \mbox{ on } (0, \infty ); \ \gamma_n (s) = (3-t) \vee 2+\frac{1}{n} \mbox{ on } (0, \infty ).
$$
Then $K'(2) \ = \ \emptyset $ but $K(2)$ consists of the function $(\gamma,a)$
$$
a(s) = s \mbox{ on } [1, \infty ); \ \gamma (s) = 2 \mbox{ on } [1, \infty ).
$$

 We now introduce the set $ \hat {K}$, the prospective limiting collection of aged paths.
For all natural numbers $n$ given $f_n (= (\gamma_n  , a_n )) \ \in \ K(n)$
 we say that $g \  \in \  K(n+1) $ is a parent of $f_n$ if $\Pi_n(g) \ = \ f_n$.  Generally $f_n$ will have multiple parents but it must have at least one.  A sequence $(f_n)_{n \ge n_0}$ is a {\it chain} if for each $n \ge n_0, \ f_n \ \in \ K(n)$ and $ f_{n+1}$ is a parent of $f_n$. It is a {\it good chain } if in addition there exists a strictly increasing integer sequence $(r_k)_{k \geq 0}$ with  $n_0 = r_0 $ and so that for each $k, \ f_{r_k} $
 has the property that $ a(s) \ \in \ [2^{-2r_{k-1}},2^{-3r_{k-1}/2}]$ 
for some $s$ in its domain of definition.  
Obviously for a  chain, the domain of the $f_n$ increases with $n$ and the $f_n$s agree on the intersection of their domains,  moreover the domain is bounded below, for example by $-n_0 - M_{n_0}$, so we can define an
$f$ on the limit of these domains $(\sigma, \infty )$ so that for each $n, \ f_n \ = \ \Pi_n(f)$. 
It is clear that for a good chain $\lim_{s \downarrow \sigma} a(s) \ = \ 0$.  To complete this part of the construction of the limit we show that
for any $n_0$ and any  $f_{n_0} (= (\gamma _{n_0} ,  a_{n_0}  )) \ \in \ K({n_0} )$, there exists a good chain $(f_n)_{n \ge n_0}$.
We then take  $ \hat {K}$ to be the collection of aged limiting paths $f$ produced by good chains, taking the union over  $n_0 \geq 1$.

To specify a good chain it is only necessary to supply a sequence $r_k$ and functions $f_{r_k}.  $  We define these recursively.  
We suppose without loss of generality that $M_t \geq  \ 2t\  \forall  \ t$.  We first choose $r_1$ via property \ref{v} to equal 
$M_{r_0}$ and $f_{r_1}$ to be the claimed function $(\gamma',a')$ for  $t=r_0$ and $(\gamma,a) \ 0 \ f_{r_0}$. 
Now we proceed recursively:
Given $r_k$ and $f_{r_k}$ with, for some $s$ in its domain \\
$a(s) \in [2^{-2r_{k-1}}, 2^{-3r_{k-1}/2}]$ . 
We simply take $r_{k+1} $ to be $M_{r_k}$.  We apply \ref{v} to function $f_{r_k}$
This gives a function $ f_{r_{k+1}} =(\gamma ',a ')$ in $K(M_{r_k} )= K(r_{k+1} )$ which has $s$ in its time domain with     $a(s) \ \in \ [2^{-2r_k}, 2^{-3r_k/2}]$ and has 

%
$\Pi_{r_k} (f_{r_{k+1}}) \ = \ (f_{r_{k}})$. 

For any $t\in[1,\infty)$, $|\Pi_t\hat{K}|\leq|\Pi_t(K(N))|\leq M_t$. by right continuity of $M$.

We similarly check that criteria \ref{ii}-\ref{v} are satisfied by $\hat{K}$.


We thus have that  $\hat{K}\in\mathcal{K}$. 
We also note that conditions \ref{i} and \ref{v} ensure that for each $t \geq 1$ the maximal functions in $\Pi_t(\hat{K})$
are indeed of the form $\Pi(f)$ for $ f \in \hat{K}$.



It remains to show that $K_j \stackrel{d_{\mathcal{X}}}{\rightarrow }  \hat {K}$.

Given the metric and bounded convergence, it is enough to show that we can find countable $V \subset [1, \infty )$ so that for $t \in V^c $, 
$$
d(K(t), \Pi_t(K_j)) \ \rightarrow \ 0.
$$
 since $\Pi_t \hat {K} = K(t)$.
We take $V = \cup _{k=1}^5 A(i) $
where $A(1) = \{ 1=t_0 < t_1 < t_2 \cdots t_m < \cdots \}$is the set of points of discontinuity for the increasing function
$t \rightarrow |K(t) |$.
For $i\geq 1$, let $m_i = |K(t) |$ for $t \in (t_{i-1}, t_i)$ and let
$f^{i}_l, \ l = 1,2 \cdots m_i$ be functions defined on $[b_{t_i, l}, \infty)$ so that for $t \in (t_{i-1},t_i)$,
$$
K(t) =\{ \Pi_t f^{i}_l\}_l
$$
\\
Similarly define $m'_k$ and $v_k$ and functions $f^k_{v_k}, \ v_k = 1,2,\cdots, m'_k$ so that
$$\Pi_t K_k = \{ \Pi_t f^k_{v_k} \}_{v_k}$$

We let $A(2) \cap (t_{i-1},t_i)= $ the union  over $l$  of points of discontinuity of $f^{i}_l $ \\
We let $A(3)$ be points $t$ in $ (t_{i-1}, t_i)$ for which $t$ ot $-t$ is a jump value for $\gamma^{i}_l$ or
$2^{-t}$ is a jump value for $a^i_{l}$ :
where $x$ is a jump value for cadlag $g$ if for some jump time $s$ for $g$, $x = g(s)$ or $g(s-)$.\\
We let $A(4)$ = the set of $t\in  (t_{i-1},t_i)$ which for some $l$ is a point of
discontinuity for $t \rightarrow b_t(f^i _l)$.\\
We let $A(5) $ be the set of $t $ so that $a^i_l(s) = 2^{-t}$ for $s$ a jump time of $\gamma ^i_l$ or
$\gamma^i_l (s) = t $ or $-t$ for $s$ a jump time for $a^i_l $ for some $l$.

For $t \in V^c \cap (t_{i-1}, t_i)$, to show that $d(K(t), \Pi_t(K_k)) \ \rightarrow \ 0$ we need only show that
for each $l, d(\Pi_tf^{i}_l,\Pi_t f^k_{v_k}) $ tends to zero where $d(\Pi_Nf^{i}_l, \Pi_N f^k_{v_k}) $ tends to zero.  
Now we have that by the definition of $A(4)$,   $s \rightarrow b_s(f^{i}_l)$, is continuous at $t$, and so if

$\Pi_tf^{i}_l,$ is continuous at $b_t(f^{i}_l)$, then we have easily that $b_t(f^k_{v_k})$ tends to $b_s(f^{i}_l)$
as $k$ tend to infinity.  In this case the desired convergence follows from Lemma \ref{lem: cadlag}.

We now suppose there is a discontinuity of $\gamma^i_l $ at $b_s(f^{i}_l)$.  In this case we must have 
$a^{i}_l(b_s(f^{i}_l) > 2^{-t}$  and continuous at this time point.  Further by the fact that $t $ is not in $A(3)$, we must have that $\gamma^i_l(  b_s(f^{i}_l) -) $ is not in $[-t,t]$ and $\gamma^i_l(  b_s(f^{i}_l) ) \in  (-t,t])$. 
Again in this case we get the desired convergence.  The case where there is a discontinuity for  $a^i_l $
follows similarly.






Hence $\hat{K}$ is the limit of subsequence $\{K_{n_j},j\in\N\}$ under metric $d_{\mathcal{X}}$.

\end{proof}

\section{Stable processes } 
We now use the general theory of aged paths to consider systems of coalescing stable processes  (for $1 <  \alpha  \leq 2$).    In \cite{MRV} a construction of a system of coalescing stable paths based on coalescing stable processes starting from dyadic rational spacetime points with  corresponding age processes was given.   
We can easily see that this system is a locally finite system of aged paths. We denote this random collection of aged paths taking value in $(\mathcal{X},d_{\mathcal{X}})$  by $\mathcal {S}$. We wish to prove that a system of coalescing random walks with the same kernel and in the domain of attraction of a 
symmetric stable process will, when suitably renormalized converge to our locally finite stable aged paths $\mathcal {S}$.

The key idea in establishing the Brownian web as the object to which scaled coalescing random walk paths from every spacetime lattice point converge is to define it as a space consisting of a collection of compact set of paths. To show the compactness of the collection one proceeds by showing the tightness of the collection of coalescing Brownian paths starting from a dense countable subset of the spacetime. This is obtained  by showing the equicontinuity of the collection. The  convergence of random walk paths to the Brownian web proceeded by   showing that 1) the limiting random set of paths  contains a subset of paths distributed as  the Brownian web almost surely and 2) the limiting random set of paths does not contain paths other than those in the Brownian web almost surely. The first follows from applying Donsker's invariance principle and the second requires some additional conditions. 

In our case,  as already noted,   we are dealing with a wilder collection of paths which can make large jumps 
so that equicontinuity is not feasible.  To deal with this we introduce an augmented path space namely the space of aged paths with a corresponding topology where convergence follows from showing convergence of the coalescing paths restricted to finite spacetime boxes increasing to $\mathbb{R}^2$ with corresponding ages having a strictly positive lower bound  decreasing to zero as the boxes increase to $\mathbb{R}^2$. The key idea is the observation, following  \cite{BG}, that in a coalescing system the number of paths at any time with a given nonzero age is locally finite.  As is similar to the  Brownian web, the stable web is considered as consisting of a collection of  aged paths satisfying the property that in any finite box there are only a finite number of paths with a given positive age. We obtain $\mathcal {S}$ as a collection of compact sets of paths by showing tightness in $(\mathcal{X}, d_{\mathcal{X}})$  We obtain this tightness by  showing that conditions \ref{i}-\ref{v}
 introduced in the previous section  hold with probability close to one. 

Since we use a topology which is a little different from the one used in  \cite{MRV} and our characterization of compact sets in $\mathcal{X}$ is different we will provide the necessary arguments to show that $\mathcal {S}$ is tight. The arguments showing that conditions  \ref{i},\ref{ii}, are satisfied with probability close to one are exactly as in \cite{MRV} Page 794 while the arguments for  \ref{iii} are similar to those in \cite{MRV} with additional ingredient given by Lemma 3.1. The main tool is that the density of paths  with age  at least $ a >0$ for some strictly positive $a$  is  
equal to $\frac{c}{a^{1/\alpha}}$ for some $c$.  To bound $\gamma$ one uses the property of the stable process that the probability of making large jumps in a given time goes to zero as the size of the time interval  goes to zero and the fact that the number of paths $\gamma$ in a box with age greater than $\epsilon >0$ is bounded in probability. A bound for the age is obtained by noting that a path with a large age has a small density and hence has a small likelihood of being in a finite spacetime box. The bound on the number of paths in a spacetime box can again be used to show that the oscillation $\omega(2^{-n}, 
\gamma^t, [-(t+1), t])$ for $t$ fixed, goes to zero in  probability as $n \to \infty$.  Probability estimates for condition  \ref{iv}  are  proved in Lemma \ref{difftimedd} of this section  and condition  \ref{v}  is a consequence of Proposition 4.8 of  \cite{MRV}  and is proven following Proposition \ref{prop32} below.  After we establish the tightness of  the collection of coalescing aged stable paths  $\mathcal{S}$ as a process taking values in $(\mathcal{X}, d_{\mathcal{X}})$ we will do the same for  scaled coalescing random walks. Once tightness is established convergence of collection of scaled aged coalescing random walk  paths to $\mathcal {S}$ proceeds exactly as in \cite{MRV}. Therefore we only give a sketch of  of that part of the proof below.

 We consider coalescing random walks  in the domain of attraction of stable law starting from  every point in $\mathbb{Z}^2$ where the time is scaled as $N^{-1}$ and space is scaled as $N^{-1/\alpha}$ taking values in $\mathcal{X}$. We denote this sequence by $\mathcal{W}_N$.  We want to show that for $t \geq 1, \ (\Pi_t(\mathcal{W}_N))$ converges in distribution to $\Pi_t(\mathcal{S})$. In the event this reduces to showing this for integer times $n$.

Given a spacetime box size $n$ we proceed by constructing discrete approximations $\mathcal{S}^{\theta,n'}$ of finite number of stable web paths starting in a larger box of dimension $n'$ in a grid of size $\theta$. A corresponding age process for these paths based on the collection  $\mathcal{S}^{\theta,n'}$ is also defined. We show that  given a margin of error $\sigma$ we can choose $n'$ and $\theta$ so that the paths in $\mathcal{S}^{\theta,n'}$ approximate the paths in $\Pi_n \mathcal{S}$  in the $d_{\mathcal{X}}$ metric with an error less than $\sigma$ with probability greater than $1-\sigma$. See \cite{MRV} for details.
Now we make use of an approximation similar to the one for $\mathcal{S}$ by using walks starting from $[-n',n']^2 \cap (\theta Z)^2$ and appropriately scaled denoted by $\mathcal{W}_N^{n',\theta}$ to approximate  $(\Pi_n(\mathcal{W}_N))$. We can show that $\mathcal{W}_N^{n',\theta}$ converges to  $\mathcal{S}^{\theta,n'}$  in distribution. Given a bounded continuous function $F$ on $\mathcal{H}$ we can control $|E(F(\Pi_n \mathcal{S} ) - E(F(\Pi_n \mathcal{S}^\theta ))| $using tightness of $\mathcal{S}$. Similarly we can control $|E(F(\Pi_n \mathcal{W}_N)) - E(F(\Pi_n \mathcal{W}_N^\theta))|$ using tightness of $\mathcal{W}_N$. This together with the convergence in distribution of $\mathcal{W}_N^{n',\theta}$ to $\mathcal{S}^{\theta,n'}$ gives us the desired result.

We fix a kernel  $p(.)$.  For  $1 < \alpha < 2$, we assume, as in \cite{MRV} that $p$ is symmetric and that
$$
\lim_{n \rightarrow \infty }n^{1+ \alpha} p(n) \ = \ C \ \in (0, \infty ).
$$
( for $\alpha =2$ we simply assume that $\sum_x p(x) \ = \ 0, \ \sum_x p(x) x^2 \ < \ \infty$   ).
We note  for $1 < \alpha < 2$ (see \cite{MRV} for details) that this condition is stronger than simply being in the domain of attraction of a symmetric $\alpha $ stable law.  It 
implies in particular that for any $x \ \in \ \mathbb{R}$ and $t \geq 0$, if $x^r/r^{1/ \alpha} \ \rightarrow \ x$ as $r$ tends to infinity then 
$$
P^{x^r} ( S_k \ \not= \ 0 \  \forall \ 0 \leq k \leq tr)
$$
 (where $P^{x^r}$ represents the probabilities for $p(.)$ random walks, $S_k$),
converges to  $P^x(X_s \not= 0 \ \forall \  0 \leq s \leq t)$ for the limit stable process.   From this point
the phrase stable process will denote this particular stable process.)
We now consider our system of renormalized walks. Let $S_r=\sum_{i=1}^rX_i$ where $(X_i)_{i\geq1}$ are i.i.d. random variables with kernel $P(X_1=x)=p(x),x\in\Z$. 
 For each $(x,r)\in\Z^2$, let $(S_m^{x,r})_{m\geq r}$ be a symmetric stable random walk starting from $(x,r)$ with  $S_r^{x,r}=x$. These stable random walks move independently until coalescence. More precisely, if there exists $r'\geq r$, $S^{r,x}_{r'}=y$ then $$\forall m\geq r', S_m^{x,r}=S_m^{y,r'}.$$
For $(S^{x,r}_m)_{m\geq r}$, the age is defined as
$$a^{x,r}(S)_{m+ 1/2}=\sup_{c}m+\frac{1}{2}-c,$$
where the supremum is over all possible  times $c \leq r$ such that $S^{y,c}_m=S^{x,r}_m$ for some $y\in\Z$. 
For other times it increases with unit rate constant on intervals $[m+1/2, m+3/2)$ and similarly on $[r,r+1/2)$. 
 The shift of $\frac{1}{2}$ is to satisfy condition (iii) in Definition \ref{def: age path}. This shift really plays no other role.   It was not present (or needed) in the previous paper as our definition of aged paths was not quite the same.

We can make $(S^{x,r}_m)_{m\geq r}$ into continuous  in time random walks by letting $S^{x,r}_t=S^{x,r}_m$ for $t\in[m,m+1)$. 

For a system of aged stable processes $(S^{x,r},a^{x,r}(S))$, for $(x,r)\in\Z^2$, define the renormalized system as
\begin{align*}
\Gamma^N=\bigg\{(x,r)\in I: &(\gamma^{x,r},a^{x,r})\text{ defined on }[r/N,\infty)\text{ by }\\
&\gamma^{x,r}_{s+\frac{r}{N}}=\frac{S^{x,r}_{r+sN}}{N^{1/\alpha}}\text{ and }a_{s+\frac{r}{N}}^{x,r}(\gamma) = \frac{a(S^{x,r}_{r+sN},r+sN)}{N}
\bigg\} ,
\end{align*}

where  the set of indices  $I=\{(x,r),a(S^{x,r})_{r+ 1/2}=1/2\}$. As  the sufixes are  indices there is no need for renormalization . 

Our objective in this section is to show

\begin{theorem} \label{thm4}
Under the aged path topology topology, $\Gamma^N$ defined above converges as $N \rightarrow \infty$ to the stable web
defined for stable processes arising out of the limits of the random walks above.
\end{theorem}

{\it Remark: } The proof shows that the distribution of aged paths for the stable web is also tight.

The proof of this is very similar to the proof of the Proposition 6.4 in \cite{MRV}, that said, for the case of $\alpha \ = \ 2$,
some of the criteria given for compact sets in that paper  do not hold. The differences
are in the main due to the difference in the compactness condition to be verified for Proposition  \ref{propcomp}
below.  We also have a similar but different topology and a slightly different $\Gamma^N$

 The chief part of the proof is showing compactness

\begin{proposition} \label{propcomp}
The system of coalescing aged random walks $\Gamma^N$ described above is  tight under the given topology.
\end{proposition}

We need simply to verify  that  the compactness conditions \ref{i}-\ref{v} of Theorem \ref{thm3} hold outside probability $\epsilon $ for fixed $\epsilon > 0$.


Condition \ref{iii} is similar to condition (iv) of \cite{MRV}.  The difference is the extension of the truncated functions domain.  The following addresses that

\begin{lemma}  \label{Jone}

Given  a cadlag function $f$ defined on interval $[a,b]$ having  for $n \geq 2$
$\omega (2^{-n}, f, [a,b]) < \frac{1}{4}$, for any $c< d  \in (a,b]$, the function
$g = f $ on $[c,d]; = f(c) $ on $[c-1,c)$  and $= f(d) $ on $[d, d+1]$ has 
$$
\omega (2^{-n}, g, [c-1,d+1])  \leq \omega (2^{-n}, f, [a,b])
$$

\end{lemma}
\begin{proof}
Given any partition $a=t_0 < t_1 < \cdots < t_m =b $ with $t_{i} \geq t_{i-1} + 2 ^{-n}  \  \forall \ 1 \leq i \leq m$, we can 
take partition $s_j $ of $ [c-1,d+1]$ by taking $s_0 = c-1,  \ s_1 = = t_k $ where $c \in \ [t_{k-1}, t_k)$
and thereafter
 taking $s_j =t_{j-1+k}$ until we arrive at interval  $[t_{l-1}, t_l)$ containing point $d$, which we replace by interval $[t_{l-1}, d+1]$ 
  It is immediate that 
$$
\sup_j Var(g,[s_j, s_{j+1})) \leq \ \sup_iVar(f,[t_i t_{i+1}))
$$
where $Var (h,I)$ is the variation of function $h$ on interval $I$.

The lemma follows.
\end{proof}
{\it Remark:} The lemma is useful in that in considering a cadlag function, $f$, with jumps, even though $\omega (2^{-n}, f, [a,b]) $ might be ``reasonable", we cannot suppose that it would be so on subdomains if the initial point of the subdomain is close and to the left of a jump of $f$.\\
In order to show condition \ref{iv} we need
 the following lemma.
\begin{lemma}  \label{difftime}
Given $x_0 , \delta, \gamma > 0, \exists \delta'' >0$ so that the probability that a stable process of rate equal to that of the given stable process or double. starting at $x: |x| \ge x_0 $ has no jump of size greater than $\delta $ within $\delta''$ of
$$
\tau_0 = \inf \{ t \ge 0 : X_t = 0\}
$$
is greater than $1 - \gamma $.
\end{lemma}

{\it Remark:} For the case  $\alpha =2$ the lemma is trivial as there are no jumps. 

\begin{proof}
 It is only necessary to treat the given stable process directly.
Suppose $(X^x_t)_{t\geq0}$ is a stable process starting from $x\in\R$   We first note that by the strong Markov property applied at time $\tau _0$ 

\[
P \left( \exists \mbox{ a jump by } X^x\mbox{ more than }\delta \mbox{ in } [ \tau_0, \tau_0 + 
\delta'' ) \right)  \le \frac{K \delta''}{\delta^ \alpha }  < \gamma / 3
\]
if $\delta '' $ was fixed sufficiently small.  It remains to deal with interval $ [ \tau_0 -
\delta '' , \tau_0 )$.  We first choose $\delta''' > 0 $ so that $\frac{K \delta'''}{\delta^ \alpha }  < \gamma / 9$ and
\[
P^y (\tau_0 > \delta ''' ) < \gamma / 6 , \ \forall | y | < \delta '''.
\]
That this can be done follows as, from scaling, $\sup_{ | y | < \delta '''}P^y (\tau_0 > \delta ''' ) \ = \ \sup_{ | y | < 1}P^y (\tau_0 > 1/ (\delta ''' )^{\alpha -1})$. 
So, as above, the probability that starting from such a $y $, that there is a jump in interval $[0 , \tau_0 ) $ is less than $\gamma / 6$.
To complete the proof, we now choose $ \delta ''$ (that satisfies the earlier inequality) so that for all $x $ with $|x| \ge x_0 $

\[
\sup _{|y| \ge \delta '''} P^y( \tau _ 0 < \tau_{\delta'''} + \delta'') < \gamma / 2
\]
for $\tau_{\delta'''} \ = \ \inf \{ t: |X_t |  \le \  \delta ''' \}$.
This can be done as given $X_{\tau_{\delta'''} -}, $, the value $X_{\tau_{\delta'''} }, $ simply corresponds to adding to $X_{\tau_{\delta'''} -}, $
a random jump so that  $X_{\tau_{\delta'''} } $ is in $[- \delta''', \delta''']$.  Give the montonicity properties of the Levy measure
for $X$, we easily have that
uniformly over $X_{\tau_{\delta'''} -}, $, P($X_{\tau_{\delta'''} } \in [- \delta_0, \delta_0]) \le \frac{2\delta_0 }{\delta''' +\delta_0}     $
for any $\delta_0$.  We choose $\delta_0$ so that this probability is less than $\gamma /9$ and then reduce $\delta ^{''}$ so that for each $y $ with $|y| \geq \delta _0,  P^y( \tau_0 < \delta^{''} ) \leq \gamma/3$.
\end{proof}
As noted, our rescaled random walks do not merely converge to $(X_t)_{t\geq0}$  in distribution, but also have hitting times of zero converging appropriately.  From this we easily obtain

\begin{lemma}  \label{difftimedd}
Given $x_0 , \delta, \gamma > 0, \exists \delta''>0$ so that the probability that for all N, the random walk $(S_n)_{n \ge 0} $ starting at $x_N$ with $x_N: |x| \ge x_0 N^{1/ \alpha }$ has no jump of size greater than $\delta N^{1/ \alpha }$ within $N \delta''$ of
$$
\tau^N_0 = \inf \{n \ge 0 : S_n= 0\}
$$
is greater than $1 - \gamma $.  The same result holds if $(S_n)_{n \ge 0} $ is replaced by  $(S_{2n})_{n \ge 0} $.
\end{lemma}

{\it Proof of Proposition \ref{propcomp}} \\

As noted, we must simply show that, given $\epsilon > 0$, constants $M_t, \  \  \delta_t (m)$ can be found so that
outside probability $\epsilon $ conditions \ref{i}-\ref{v} hold.  Which will be shown if we can find 
for each condition a suitable choice so that the condition holds for each $\Gamma^N$ outside probability $\epsilon / 5$.
   The conditions \ref{i} and \ref{ii} are exactly as in \cite{MRV}, page 794 (i)-(iii).   
We begin with condition \ref{i}.  By monotonicity of $|\Pi_t(K)|$ in $t$ it is sufficient to show that we can find $M^1_n, \ n \geq 1$ so that
for each $N$
$$
P(|\Pi_n(\Gamma^N)| \ > \ M^1_n) \leq \ 2^{-n}\epsilon / 5.
$$
We have by the classical (and very robust) \cite{BG} argument that the coalescing system $\left( (S^{x,0}_r)_{r \geq 0} \right)_{x \in \mathbb{Z}}$
is such that the random set $\{ (S^{x,0}_r)\}_{x \in \mathbb{Z}}$ has density bounded by $c/r^{1/\alpha}$ for fixed $c$  not depending on $r$. So
for each $s$ and each $N$, $\{ \frac{ ( S^{x,0}_{sN})}{N^{1/ \alpha}}\}_{x \in \mathbb{Z} }$ has density bounded by
$\frac{c}{(sN)^{1/\alpha}} N^{1/ \alpha} \ = \ c/s^{1/\alpha}$.

We consider for $-n2^{(n+2)}-1 \leq i \leq n2^{(n+2)}$, the random set
$$
V^{N,n,i} \ = \  \{ \frac{(S^{x,iN2^{-(n+2)}}_{(i+1)N2^{-(n+2)}})}{N^{1/ \alpha}}\}_{x \in \mathbb{Z}} \cap [-n,n]
$$
We have immediately that for the above $c$, $E|V^{N,n,i} | \ \leq \ c2^{(n+2)/\alpha }(2n+1)$
and so $\sum_{i=-n2^{(n+2)}-1}^{n2^{(n+2)}} E|V^{N,n,i} | \ \leq \ c2^{(n+2)/\alpha }(2n+1) 2(n2^{(n+2)}+1)$. for each $N$.

This bound does not quite give a bound for the expectation of $|\Pi_n(\Gamma^N)|$: it is possible for a $(\gamma,a) \in \Gamma^N$ to have 
$\Pi_n((\gamma,a))$ non null but such that for every time $i2^{-(n+2)}$ for $i$ as above, $ \gamma (i2^{-(n+2)})$  (when defined) not be in $[-n,n]$.

However (for $(\gamma,a) \in \ \Gamma^N$) if $\Pi_n((\gamma,a))$ is non null and if
$$
b_n \ = \ \inf \{ s: \gamma_s \in [-n,n], \ a_s \geq 2^{-n} \},
$$
 then for
$$
i \ = \  \lfloor (b_n - 2 ^{-(n+1)}) 2 ^{(n+1)} \rfloor ,
$$
which is in $[-n-2 ^{-(n+1)},n]$,
there must exist $x \in \mathbb{Z}$ so that 
$$
 N^{1/\alpha} \gamma(b_n+s) \ = \ X^{x,(2i -1)2 ^{-(n+2)} N}_{(b_n+s)N} \forall s \geq 0.
$$
From this  we see that $E(\vert \Pi_n(\Gamma^N \vert )) \ \leq \  \sum_{i= -n2^{(n+2)}-1    }^{n2^{(n+2)}} 
E(|V^{N,n,i}_-|)$, where $V^{N,n,i}_-$ is the set of the
$ \{ \frac{(S^{x,(2i-1)N2^{-(n+2)}}_{(i) N 2^{-(n+2)}})}{N^{1/ \alpha}}\}_{x \in \mathbb{Z}}$
so that for some $s \ \in \ [iN2^{-(n+1)} ,iN2^{-(n+1)} +N2^{-(n+1)}]$,
$S^{x,(2i-1)N2^{-(n+2)}}_s /N^{1/ \alpha} \ \in \ [-n,n]$.  Let 
$$
c_1 \ = \ \inf_{N,n} \  \inf_{x\in [-nN^{1/ \alpha}, n N^{1/ \alpha}]} \inf_{s\in [0,2nN]} P`^x(S_s /N
^{1/ \alpha }\in [-n,n])
$$
which is easily seen to be strictly positive.  
For our system $\Gamma ^N$ consider the natural  filtration.

We have by the Markov property applied to this filtration at the stopping times when walks hit $[-n,n]$ that 
$$
E(|V^{N,n,i}_-|) \ \leq \ \frac{1}{c_1}E(|V^{N,n,i}|) \ \leq \ \frac{1}{c_1}c2^{(n+2)/\alpha }(2n+1).
$$
So uniformly in $N, E(|\Pi_n (\Gamma^N)|)  \ \leq \ \frac{1}{c_1}c2^{(n+2)/\alpha }(2n+1) 2(n2^{(n+1)}+1)$.

We may therefore take $M^1_n$ to be $\frac{ 2^n 5}{\epsilon}\frac{1}{c_1}c2^{(n+2)/\alpha }(2n+1) 2(n2^{(n+1)}+1)$.

For property \ref{ii},  there are two probabilities to bound for two different events relating to $\gamma $ and to $a$.   So we need to bound the probabilty of each by $\frac{\epsilon}{2^n 10}$.   We first first consider the spatial variation of paths $(\gamma, a)$ with $\Pi_t((\gamma, a))$ non null.  
Again by montonicity, it is enough to show that there exists $M^{2,a}_n$ for integers $n \geq 1$ so that  for each $N$,
$$
P(\vert \Pi_n(\Gamma ^N) \vert\leq M^1_n, \exists (\gamma, a) \ \in \   \Pi_n(\Gamma ^N) : \vert \gamma _s \vert > M^{2,a}_n \mbox{ for }
s \in [b,n]) \ \leq \ \epsilon/( 2^n 10).
$$
By weak convergence, we have that for every $\delta > 0, \ \exists L \ = \  \ L(\delta ) $ so that
$$
P( \sup_{0 \leq  s \leq 2nN}  \vert S^{0,0}_s \vert \geq LN^{1/ \alpha}) \ < \ \delta.
$$

For our system $\Gamma ^N$ consider the filtration $\mathcal{F}_j \ = \  \sigma \{S^{x,k}_r : x \in \mathbb{Z}, k \leq r \leq j\}$.
We  use the Markov property of our random walks applied  at time $b$ for new walks with  $ \Pi_n(\gamma, a)$ non null 

$$
P(\vert \Pi_n(\Gamma ^N) \vert\leq M^1_n, \exists (\gamma, a) \ \in \   \Pi_n(\Gamma ^N) : \vert \gamma _s \vert >(L+n)N^{1/ \alpha  } \mbox{ for }
s \in [b,n]) \ \leq \ M^1_n \delta.
$$
Accordingly we choose $M^{2,a}_n \ = \ L(\delta )+n$ for $\delta < \epsilon /(M^1_n 2^n 10)$.

To control the ages, we need to find $M^{2,b}_n$ so that for all $N$
$$
P(\vert \Pi_n(\Gamma ^N) \vert\leq M^1_n, \exists (\gamma, a) \ \in \   \Pi_n(\Gamma ^N) : \vert \gamma _s \vert \leq M^{2,a}_n \mbox{ for }
s \in [b,n], a_s \geq M^{2,b}_n ) \ \leq \ \epsilon/(2^n 10 ).
$$

We again use the uniform in $N$ bound; $\{ (\frac{S^{x,0}_{sN}}{N^{1/ \alpha}})\}_{x \in \mathbb{Z} }$ has density less than or equal to $c/s^{1/\alpha}$.
Thus the density of points $\{  (\frac{S^{x,-M^{2,b}_n N +nN}_{nN}}{nN^{1/ \alpha}}) \}_{x \in \mathbb{Z}}$ is less than $\frac{c}{(M^{2,b})^{1/ \alpha}}$.
In particular, the probability that this set intersects $[-M^{2,a},M^{2,a}]$ is bounded uniformly in $N$ by
$$
\frac{(2M^{2,a}+1)c}{(M^{2,b})^{1/ \alpha}}.
$$
Again, this does not quite give the required bound on the probability that an aged path can intersect spacetime rectangle $[-n,n] \times [-n,n]$.  However as before we have that the probability that  there exists $x, \ \frac{S^{x,-M^{2,b}_n N +nN}_{.}}{nN^{1/ \alpha}}$intersects $[-n,n]$ for time in $[-nN,nN]$ is bounded by $ \frac{(2M^{2,a}+1)c}{c_1(M^{2,b})^{1/ \alpha}}. $
This is a bound for the probability of the event $\{ \vert \Pi_n(\Gamma ^N) \vert\leq M^1_n, \exists (\gamma, a) \ \in \   \Pi_n(\Gamma ^N) : \vert \gamma _s \vert \leq M^{2,a}_n \mbox{ for }
s \in [b,n], a_s \geq M^{2,b}_n \}$ and we simply take $M^{2,b}_N$ so large that it is less than $\epsilon/(2^n 10)$.
Define $M^2_n := M_n^{2,a} \vee M_n^{2,b}.$ 

Condition \ref{v} follows from a version of
Proposition 4.8 of \cite{MRV}:
\begin{proposition} \label{prop32}

Fix $1 \leq n \ < \ \infty $ and $\delta > 0$.
Given $\eta_2,\ \sigma  >0$,  there  is $M \ = \ M(n,\delta, \eta_2, \sigma) $ so  that  the  probability  that  a  path $(b,\gamma,a)$ (for $\Gamma$ the collection of aged paths that is either the stable web or $\Gamma^N$ the renormalized random walks) satisfies \\
(i)$ \gamma $ hits $[-n,n]^2 $ while of age at least $ \delta$,\\
(ii)  there exists $t \in $ $ [-n,n]$ so that $ a(t) \geq \eta_2 $ and $ \gamma (t) \in \ [ -M,M]^c $ \\
is bounded  above by $\sigma $.

\end{proposition}

{\it Remark:} Proposition 4.8 in \cite{MRV} explicitly has $\delta $ equal to $2^{-n}$ but this specific choice of $\delta $ plays no key role in the argument
for the existence of $M$.

To show \ref{v} it is only necessary to supply an $M_5$ for a given positive integer $n$ and not for continuous $t \geq 1$ by monotonicity as before. Given $n$ we take $\delta $ to be $2^{-n}$ and $\sigma $ to be $\epsilon / 5$ in our application of the above proposition but 
$n$ is replaced by $n' = n+M^2_N$ above.   $\eta_2 $ will be $2^{-2n}$.  If an aged path $(\gamma,a)$ has $\Pi_n(\gamma,a) $ is nontrivial, then certainly there exists $s \in [-n,n]$ so that $\gamma(s) \in [-n',n']$ and $a(s) \geq \epsilon \ = \ 2^{-n}$.  Necessarily its domain of definition is contained in $(-n-M^2_n, \infty ) \ = \  (-n', \infty )$ by property (i) above.  Applying the   proposition gives the result for the first part of \ref{v}.  For the second part we simply note that all the choices of the collection of aged paths give the existence of an aged path
$(\gamma', a')$ so that $\Pi_n(\gamma', a') \ = \ \Pi_n(\gamma, a)$ and such that on $\sigma' < s \leq \sigma'+2^{-n}, \ a(s) = s - \sigma' $.
This immediately yields the second part.

For condition \ref{iii}, given $t \geq 1$ and  $g \in \Pi_{t} \Gamma^m $ with domain increased as in \ref{iii}, let
$n$ be an integer greater than $t+2$.
We apply Lemma \ref{Jone} to  $f \in \Pi_{n+2} \ \Gamma^m $ corresponding to  $g  $ with $[a,b]
\ = \ [b_{n+2}, n+2]$,  $c= b_t$ and $d= t$ to obtain
$$
\omega (2^{-n}, g, [-t-1,t])  \leq \omega (2^{-n}, f, [-n,n]) = \delta'_t(n)
$$

\ref{iii} now follows from the proof of condition (iv) of page 794 of \cite{MRV}.

The main work is addressing condition \ref{iv}. Again note that for $\alpha = 2$, there is nothing to be shown.

We note that the age process can only jump via a coalescence.  Given $n$, it is sufficient to find $\delta_n$
so that the probability that there exists path $(b_n,\gamma,a)$ so that
\begin{itemize} 
\item
$\Pi_n(\gamma,a)$ is not null 
\item
$\gamma$
jumps by at least $2^{-n}$ within time $\delta_n$
of  a coalescence with another aged path of age at least $2^{-n}$ (at the time of coalescence)

\end{itemize} 
is smaller than $\epsilon / 2^n 10$ for all  $n$ 
(uniformly in $N$).
The complement of this event is contained in the union
$$
\cup_{-n2^{n+2}}^{n2^{n+2}}A(i)
$$
where $A(i)$ is the event that two paths that touch $[-n,n]^2$ coalesce at time in 
$[\frac{i+1}{2^{n+2}}, \frac{i+2}{2^{n+2}}]$ and both have age at least $2^{-(n+1)}$ at time 
$\frac{i}{2^{n+2}}$.
Arguing as in steps proving B and E, we can find 
first $M^{2'}_n$, then
 $M_n ^ 3 $ and then $M_n^4 $ so that outside probability $\epsilon /( 2^n 10 ) $
\begin{itemize} 
\item
No $(\gamma,a)$ with nontrivial $\Pi_n$ projection has  $\vert \gamma(s) \vert  \geq M^{2'}_n$ for 
$s \in [-n,n]$
\item
No 
$(\gamma,a) $ which hits $[-M^{2'}_n,M^{2'}_n]$ with age at least $2^{-n}$ has age at least
$2^{-(n+1)}$ while $  \gamma \in [-M^3_n,M^3_n]^c$, $s \in [-n,n]$
\item
The number of paths which enter $[-M^3_n,M^3_n]$ while having age at least $2^{-(n+1)}$ 
is bounded by $M^4_n$. 
\end{itemize} 

We now abuse notation by replacing events $A(i)$ by the intersection with the former event $A(i)$ with the 
intersection of the three events above.  We pick $x_0$ so that a stable process (or renormalized random walk) beginning with initial value $x \in [-x_0,x_0]$ will hit zero before time $2^{-(n+2)}$ outside probability
$$
\frac{\epsilon }{4( 2^n 10 ) 2n2^{n+2} (M^4_n)^2}
$$
Given this $x_0$, we choose $\delta_n $ be provided by Lemma \ref{difftimedd} ( applied to the process that is the difference of two independent walks )with values $x_0, \delta = 2^{-n}
$ and $\gamma = \frac{\epsilon }{4( 2^n 10 ) 2n2^{n+2} (M^4_n)^2}$ and $\delta'' = \delta_n$.  (We can also suppose that $2^{-n-2} >\delta_n $).

Then we easily see that 
$$
P(A(i)) \ \leq \ \frac{\epsilon }{2( 2^n 10 ) 2n2^{n+2}   }
$$
We let $\delta_t(n) = \delta'_t(n) \vee \delta_n$ and
we are done.



 {\it Proof of Theorem \ref{thm3}}
To show show convergence of systems of coalescing random walks in the domain
of attraction for a symmetric stable process to the corresponding stable locally finite collection of aged paths, having established tightness,
we can now proceed as in \cite{MRV} leading to Proposition 5.3 as sketched out earlier.

\section{Brownian case}
In this section we restrict attention to the Brownian stable web
and establish two path properties for the brownian stable web. 
We believe that they hold in greater generality but can only prove them in the brownian case due to 
being able to exploit the path continuity and Gaussian tails.
\begin{proposition}
\label{prop: brown 1}
For the Brownian stable web, each aged path $(\gamma,a)$ defined on $(\sigma,\infty)$ has the property that $\lim_{s\downarrow\sigma}\gamma(s)$ exists. That  is we can extend $(\gamma,a)$ continuously to $[\sigma,\infty)$.
\end{proposition}

\begin{proposition}
\label{prop: brown 2}
For the Brownian stable web extended by Proposition \ref{prop: brown 1}, 

with probability one there does not exist $t \in \mathbb{R}$ and aged paths $(\gamma,a)$ and  $(\gamma^\prime,a^\prime )$ so that
\begin{itemize}
\item
$t > \sigma, \sigma ^\prime$, the respective birth times,
\item
$\gamma(t) \ \ne \  \gamma^ \prime (t)$
\item
$ a(t) \ =  \ a ^\prime (t)$.
\end{itemize}

\end{proposition}
 In discussing the age of a path, we exploit the coalescing property of our model and the fact that when two paths meet the age of the ``older" path is unchanged, while that of the ``younger" path jumps to that of the older.  The effect of this is that as agents affecting the evolution of the ages of paths, an aged path is effectively killed upon hitting a (strictly) older path.  We have in consequence that the age of a path at a given time, $t$, will correspond to a path whose age has never jumped.  That is to a path which has birth time $\sigma < t$ and whose age at $ s \in (\sigma, t]$ will equal $s- \sigma$.  Though this section deals with 
the special case $\alpha = 2$, this property is so for all the coalescing models.

The proofs of both these results depend on the path continuity of brownian paths, which entails (as with the brownien web but unlike for the stable web with $\alpha < 2$) that all paths that pass through a point in a compact
space time region must be within some (typically larger) compact space time region over any bounded time interval.  

Mathematically this is captured in the following ``building block" result.

\begin{lemma}
\label{lem: brown 1}
Let $A_N$ be the event that  for each
 $2\leq k\leq N$
there does not exist a path $(\gamma,a)\in\mathcal{S}$ so that $\sigma< -1$,
$$\gamma(0)\in (k-1,k), a(0)>1\text{ and }\gamma(s)\in (k-1,k),\forall 0\leq s\leq 1$$.
Then there exist constants $c,C\in(0,\infty)$ so that $$\P(A_N)\leq Ce^{-cN}.$$
\end{lemma}
{ \it Remark: } The requirement $a(0) >1$ is arbitrary and any strictly positive value will do in practice.
\begin{proof}
Consider the event $V_k$ that the Brownian path $B^{k-1/2,-1}_{\cdot}$ starting at $(k-1/2,-1)$ has the property that it does not remain in $(k-1,k)$ for $t\in[-1,1]$. Then $$\left(\bigcap_{k=2}^NV_k\right)^c\subset A^c_N.$$
But as $V_k$ can be written in terms of  independent  Brownian motions killed upon leaving disjoint intervals, we get
$$\P(A_N)\leq\prod_{k=2}^N\P(V_k)=(1-c)^{N-1}$$
 for some $0<c<1$. 
\end{proof}
{\it Remark:}
By reflection symmetry,
$$P(A'_N)\leq Ce^{-cN}$$
for $A'_N$ the event that for all  $2\leq k\leq N$  there is no Brownian path in $(-k,-k+1)$ at time $-1$ that remains in $(-k,-k+1)$ in time interval $[-1,1]$.

{\it Remark:}
By scaling this result implies that on any bounded space time region enclosed in $[a,b] \times   [s,t]  $ 
with high probability there will be "insulating" paths that are to the left of $a$ (respectively to the right of $b$) on open time intervals containing
$[s,t]$ on which they are defined.

\begin{proof}
[Proof of Proposition \ref{prop: brown 1}]

It suffices to show that for all closed finite intervals $I_1, \ I_2$ and $I_3$ with $I_1$ and $I_2 \ \subset \ \mathbb{R}$ 
and $I_3 \ \subset \ (0, \infty )$, $\forall \epsilon>0$, with probability greater than $1-\delta(\epsilon)$, if
$$
(\gamma,a) \in\Gamma_1   \equiv   \left\{ (\gamma',a')\in \mathcal{S}:(\gamma'(s),a'(s))\in I_1 \times I_3 \text{ for some } s \in I_2 \right\},
$$
then $\sup_{\sigma<s,s'<\sigma+\epsilon}|\gamma(s)-\gamma(s')|\leq \eta(\epsilon)$ for some $(\eta(\epsilon), \ \delta(  \epsilon ))\downarrow0$ as  $\epsilon\downarrow 0$.  The result then follows as all paths belong to $\Gamma_1$
for some $I_i$ in the countable collection of intervals with rational endpoints.
For definiteness we choose $I_i \ = \ [1/2,1]$ for $i = 1,2 ,3$ but it will be clear that the argument extends to any such rectangles.

It is only necessary to treat $\epsilon$ of the form $2^{-n}$ for $n\in\N$.
The following is a consequence of 
the definition of  insulating paths and path continuity. 

\begin{corollary}
Let $A_N, A'_N $ be the events of Lemma \ref{lem: brown 1} and the succeeding remark.
On event $(A_N\cup A'_N)^c$, every aged path $(\gamma,a)\in \Gamma_1$ satisfies $\gamma(s)\in [-N,N],\forall s\in (\sigma,1]$.
\end{corollary}

Given an aged path $(\gamma,a)$ denote $\lambda\left(2^{-n}\right)=\inf\{s>\sigma:a(s)\geq 2^{-n}\}$. Note that by (ii) of Definition 2.1,
$$\lambda\left(2^{-n+1}\right)-\lambda\left(2^{-n}\right)\leq 2^{-n}.$$
If $(\gamma,a)\in\Gamma_1$, then there exists $i\in\N$ so that $i2^{-(n+2)}\in [- 1/2,1]$ and $0<\lambda(2^{-n})-i2^{-(n+2)}\leq 2^{-(n+2)}$. Hence
$$\left[\lambda\left(2^{-n}\right),\lambda\left(2^{-n+1}\right)\right]\subset\left[i2^{-(n+2)},(i+1)2^{-(n+2)}+2^{-n}\right] ,$$
since $\lambda(2^{-n+1}) \leq (i+1) 2^{-(n+2)} + 2^{-n}$. 
On $\left[\lambda(2^{-n}),\lambda(2^{-n+1})\right]$, by maximality $\gamma$ must equal $\gamma^{'}$ for some $\gamma^{'}$ whose age at $i2^{-(n+2)}$ is at least $2^{-n}-2^{-(n+2 )}\geq 2^{-n-1}$
and  on event $(A_N\cup A'_N)^c$,    $\gamma^{'}(t)\in [-N,N]$ for  $t \in(-1/2,1]$. So we have that the event
$$
B(n)=\left\{\exists(\gamma,a)\in\Gamma_1,\sup_{s,s'\in[\lambda(2^{-n}),\lambda(2^{-n+1})]}|\gamma(s)-\gamma(s')|\geq n2^{-n/2}\right\}\bigcap (A_N\cup A'_N)^c
$$
$\subset \bigcup_{i2^{-(n+2)}\in[-1/2,1]}A(i,n)$,
where
\begin{align*}
A(i,n)=\bigg\{&\text{At time }i2^{-(n+2)}, \text{ a path of age at least }2^{-(n+1)}\text{ satisfies }\\
&\sup_{s,s'\in\left[i2^{-(n+2)},i2^{-(n+2)}+2^{-n+1}\right]}|\gamma(s)-\gamma(s')|\geq n2^{-n/2},\gamma(i2^{-(n+2)})\in [-N,N] \bigg\}.\end{align*}

Brownian tail probability and the expectation for the number of noncoalesced paths in $[-N,N]$ 
at time $i2^{-(n+2)}$ (See Proposition 2 of \cite{MRV} ) with age at least $2^{-(n+1)}$ give that
$$\P(A(i,n))\leq CN2^{n/2}e^{-cn^2}.$$
So
\begin{align*}
\P(B(n))&\leq 6CN2^{3n/2}e^{-cn^2}\\
&\leq C'Ne^{-cn^2/2}
\end{align*}
for universal constants $c,C'>0$. Thus
\begin{align*}
&~~~~\P\left(\exists (\gamma,a)\in\Gamma_1\text{ s.t. }\sup_{\sigma<s,s'<\sigma+\lambda (2^{-n_0})}|\gamma(s)-\gamma(s')|\geq\sum_{n\geq n_0}n2^{-n/2}\right)\\
&\leq\sum_{n\geq n_0}NC'e^{-cn^2/2}+\P(A_N)+\P(A'_N)\\
&\leq e^{-cn_0}.
\end{align*}
This bound (and maximality) implies the claim.
\end{proof}

We now turn to the proof of Proposition \ref{prop: brown 2}.
The proof comprises two parts.  We need the following definition:\\
\begin{definition}
A point $(x,t)$  gives birth to an $\epsilon $-path (for $\epsilon > 0$) if there is an aged path  $( \gamma,a) \in \mathcal {S}$ with
$\sigma \ = \ t$, age equal to $s-t$ on $(t,t+ \epsilon)$ and $\lim_{s \downarrow t}   \gamma (s) \ = \ x$.  
In this case we say that (whatever the value of $ \epsilon $) that $(x,t)$ gives birth to $(\gamma,a)$: A point
$(x,t)$ gives birth to a path  if it gives birth to an $\epsilon $-path for some $\epsilon >0$.
\end{definition}

The discussion following the statement of Proposition \ref{prop: brown 2} entails that
for any path $(\gamma^{'}, a^{'})$ (and $t > \sigma^{'} $ with (extended domain $[\sigma^{'}, \infty)$), there must be some  $(\gamma,a)$ so that $a(s) \ = \ s- \sigma $ on $[\sigma, t ]$. and $\gamma (t) = \gamma ^{'}(t)$
and $a(s) = a^{'} (s$ This property is an important reason why there can not be two paths starting at a random space time point  $(x,t)$ as stated in Proposition \ref{prop: brown 3} below.  In the topology of Brownian web there are no such constraints on paths starting at a point which makes it possible to have two (and even three)  paths starting from a random space time point $(x,t)$.

We must first rule out the possibility that  some random space time point $(x,t)$ gives birth to two distinct paths.  Then we must show that there cannot exist random times $t$ and distinct random points $x$ and $y$ so that  both $(x,t)$
and $(y,t)$ give birth to paths.

\begin{proposition}
\label{prop: brown 3}
For the Brownian stable web extended by Proposition \ref{prop: brown 1}, with probability one there does not exist
a random $(x,t)$ at which two distinct paths are born.
\end{proposition}

\begin{proposition}
\label{prop: brown 4}
For the Brownian stable web extended by Proposition \ref{prop: brown 1}, a.s. there do not exist two paths born at the same time but at different spatial  points.
\end{proposition}

In showing the impossibility of multiple paths  being born at a single space time point,(in Proposition  \ref{prop: brown 3}) we  need the following lemmas the first of which follows easily from \cite{BG}
\begin{lemma} \label{beegee}
There exist nontrivial constants $c $ and $K_c$ so that
for interval $I \ = \ [a,b], $ and integer $m$, the system of coalescing brownian motions killed on leaving $I$, starting at full occupency at time $0$ has less than $m$ distinct points at time $(b-a)^2/(cm^2)$ outside of probabiliy $K_ce^{-cm}$. 

\end{lemma} 

We will also need.  
\begin{lemma} \label{4B}
For independent brownian motions $B^i , \ i \ = \ 1,2,3 $ and $4$ each beginning on the interval $[0, \delta]$,  the  probability that none of the motions touches another before time $1$ is less than $K \delta^{c_1}$ for some universal $K$ and $c_1 > 3$ uniformly over all initial positions in the given interval .
\end{lemma} 

Proofs of these lemmas will be given after the proofs of the propositions \ref{prop: brown 3} and \ref{prop: brown 4}. In what follows by paths we mean aged paths in $\mathcal{S}$.


Let $A$ be the event $\{   \exists \  (x,t): \mbox{  Two distinct paths are born  at }  (x,t) \}$.
It is clear that 
$$
A \ = \ \cup_r \cup_{m_1}\cup_{m_2} A(\frac{1}{r}, [m_1,m_1+1],  [m_2,m_2+1])
$$

where $m_1, m_2$ range over $\Z$ and $r$ takes values in the natural numbers
and 

where $A(\epsilon , [m_1,m_1+1],  [m_2,m_2+1])$ is the event $ \{ \exists   (x,t) \in [m_1,m_1+1] \times  [m_2,m_2+1] : \mbox{  two distinct } \epsilon  \mbox{ paths are born  at }  (x,t) \}$.


\begin{proof}[Proof of Proposition \ref{prop: brown 3}]
 Using translational invariance it is only necessary to show that $P(A(\epsilon,[0,1],[0,1]))=0$ for arbitrary $\epsilon >0$. Since
$$
A(\epsilon, [0,1],[0,1]) \ \subset \ \cup _{n_1 =0}^{2^n}\cup _{n_2 =0}^{2^{2n}} C(\epsilon,n_1,n_2)
$$
where $C(\epsilon,n_1,n_2) \ = \ $
$$
  \exists \ \{ (x,t) \in [\frac{n_1}{2^n}, \frac{n_1 +1}{2^n}] \times  [\frac{n_2}{2^{2n}}, \frac{n_2 +1}{2^{2n}}] : \mbox{  two distinct } \epsilon  \mbox{ paths are born  at }  (x,t) \},
$$
it will suffice to show that $2^{3n}  \sup_{n_1,n_2}P(C(\epsilon,n_1,n_2) ) $ tends to zero as $n$ becomes large,
which again by translation invariance is the same as showing $2^{3n} P(C(\epsilon,n_1,n_2) ) $ tends to zero for any choice of $n_1$ and $n_2$.

To avoid notational encumbrance we prove this explicitly for $n_1=n_2 = 0$.  We proceed by upper bounding  the probability of an event which contains our event. We consider a spatial interval  containing $[\frac{0}{2^n},\frac{1}{2^n}]$  which is long enough so that with high probability we can find insulating paths on the left and right of the interval $[\frac{0}{2^n},\frac{1}{2^n}]$. These insulating paths are Brownian motions with age at least $2^{-2n}$ at time $0$  which remain in a spatial interval of the form $(i 2^{-n}, (i+1) 2^{-n}) $ on the right and $(-(j+1) 2^{-n},-j 2^{-n})$ on the left during the time interval $(0,2(2^{-2n}))$.  Using this we can control the possibility of coalescence with paths of age greater than or equal to zero at time zero for the paths in the interval $(-(j+1)2^{-n},(i+1) 2^{-n})$  which started at time $2^{-2n}$ during the time interval $[2^{-2n},2 \ 2^{-2n}]$. This reduces the problem to  estimating the probability that two paths which started in $(-(j+1)2^{-n},(i+1) 2^{-n})$ at time $2^{-2n}$ do not meet each other and the two insulating paths until time $\epsilon$.  This is because the union of such events over all pairs of such paths in $(-(j+1) 2^{-n}, (i+1) 2^{-n} ) $ contains our event.

We say that $1 \leq i \leq Kn:$ is  {\it left good} if $\exists (\gamma,a)$ so that $\sigma <-2^{-2n}$ and $\gamma (s) \ \in (\frac{-i}{2^n}, \frac{-i+1}{2^n}) \ \forall 0 \leq s \leq \frac{2}{2^{2n}}$.  Here $K$ is a large fixed constant not depending on $n$.
Let event $V_1$ be the event that all $1 \leq i \leq Kn$ are not left good.   Similarly 
$1 \leq i \leq Kn:$ is  {\it right good} if $\exists (\gamma,a)$ so that $\sigma <-2^{-2n} $ and $\gamma (s) \ \in (\frac{i}{2^n}, \frac{i+1}{2^n}) \ \forall 0 \leq s \leq \frac{2}{2^{2n}}$ and  let  $V_2$ be the event that all the $1 \leq i \leq Kn$ are not right good.  
The final ``bad" event is $B_1$: the number of distinct coalescing brownian motions produced by $[-Kn \frac{1}{2^n},Kn \frac{1}{2^n} ]  \times \{\frac{1}{2^{2n}}\}$ at time $\frac{2}{2^{2n}}$, killed upon leaving spatial  interval $[-Kn \frac{1}{2^n},Kn \frac{1}{2^n} ]$ ,  is at least $Cn$ for $C$ a large constant .  How large $K$ and $C$
need to be will to be shortly specified.

From Lemma \ref{lem: brown 1} (and scaling) it follows that $P(V_1 \cup V_2) \ \leq \ e^{-4n}$ for  $n$ large if value $K$ is fixed sufficiently large.  Equally we apply Lemma 
\ref{beegee} with $a = -nK2^{-n}, \ b \ = \ nK2^{-n}$  and $m= 2Kn/c$ for $c$ sufficiently small to see that outside probability $K_ce^{-2Kn }$ by time $n^2 2^{-2n} /(n^2)$, there will be less than $2Kn/c $ brownian motions for the coalescing system on $[a,b]$ at time $2 \ 2^{-2n}$ having started at time $2^{-2n}$. We  take, as we may, $K$ to be greater than $2$.

On $(V_1 \cup V_2 \cup B_1)^c$ there are less than $2Kn/c$ brownian motions at time $ 2 \ 2^{-2n} $ that were produced .by $[a,b]$ at time $ 2^{-2n}$ and have not met a brownian path of age at least $0$ at time $0$. We enumeraute these motions 
by $B^1, B^2, \cdots B^{2Kn/c}$ in order.   In actuality on event $V_N$ the number of Brownian motions is less than  $2Kn/c$, not equal to it.  We adopt the convention that the ``surplûs" Brownian motions are not counted.   We let $i_1 $ be  the first $i$  that is left good and $i_2$ be the first $i$ that is right-good.  
Let $D(i_1,i_2, j,k)$ (for $1 \leq i_1, i_2 \leq 2Kn$  and $1 \leq j < k \leq 2 Kn/c$) be the event
\begin{itemize}
\item
$\emptyset $ if either (at time  $ 2^{-2n}) , \ j$ is to the left of a brownian motion in the $i_1$'th left interval of age at least $ 22^{-2n}$ or if
$k$ is
to the right of such a motion in the $i_2$'th right interval.
\item
otherwise it is the event that the $4$ brownian motions appearing above, $ B^j,\ B^k$ and the brownian motions corresponding to $i_1$ and $i_2$ avoid each other for the next $\epsilon - 22^{-2n}$
time units.
\end{itemize}
We have that  $P(C(\epsilon,0,0 )  \leq \sum P( D(i_1,i_2, j,k)) \ + \ K_ce^{-2Kn} + 2^{-4n} $.  By Lemma \ref{4B} (with $c ^\prime  = c_1 -3$) this gives the bound
$(4 \sqrt{ n }/ \epsilon) ^{(3+ c ^\prime)/2}4K^2 n^2 4K^2n^2/c^2 2^{-c^\prime n}  \ + K_ce^{-(2K-3)n}+ 2^{-n}  ) 2^{-3n} $ for  $P(C(\epsilon,0,0 ))$.  As noted this suffices.
\end{proof}

\begin{proof}[Proof of Proposition \ref{prop: brown 4}]
If the event
$$\left\{\exists t\in\R, (\gamma,a,\sigma),(\gamma',a',\sigma')\in\mathcal{S}\text{ s.t. }\gamma(t)\neq\gamma'(t),a(t)=a'(t)>0\right\}$$
has strictly positive probability, then given Proposition \ref{prop: brown 1},
\begin{align*}
\bigg\{\exists t>0, &(\gamma,a,\sigma),(\gamma',a',\sigma')\in\mathcal{S}\text{ s.t. }\\
&\sigma=\sigma',a(s)=s-\sigma=a'(s)=s-\sigma'\text{ for }s\in (\sigma,\sigma+t),\gamma(\sigma+t)\neq\gamma'(\sigma+t)\bigg\}
\end{align*}
has strictly positive probability. Denote  by $A^u_{a_1,b_1,a_2,b_2}(s,t)$ the event that there exist $(\gamma,a,\sigma),(\gamma',a',\sigma')\in\mathcal{S}$  satisfying the following properties
\begin{enumerate}[label=(\roman*)]
\item $\sigma=\sigma'\in(s,t)$ and $\gamma(\sigma)\in [a_1,b_1],\gamma'(\sigma')\in[a_2,b_2]$
\item $a(v)=a'(v)=v-\sigma$ on $(\sigma,\sigma+u)$
\item $\gamma(\sigma+u)\neq\gamma'(\sigma+u)$.
\end{enumerate}
If the above two events have positive probability, then for some disjoint rectangles  $[a_1,b_1]\times[s,t]$ and $[a_2,b_2]\times[s,t]$, for $t-s \leq 1$,  $\P\left(A^u_{a_1,b_1,a_2,b_2}(s,t)\right) \geq c > 0$ for $t-s \leq 1$.
This implies that $\forall N\in\N$, $\exists (s_N,t_N)\subset[s,t]$ so that $t_N - s_N \ = \ 1/N$ and
$$\P\left(A^u_{a_1,b_1,a_2,b_2}(s_N,t_N)\right)\geq\frac{c}{N}.$$
We may assume by spatial translation invariance that
$$[a_1,b_1]=[-B,-A]\text{ and }[a_2,b_2]=[A,B],$$
where $0<A<B$. We then  construct $3$ realizations of the Brownian stable web.
 The first is based on coalescing system of Brownian motions  $(B^{x_i,t_i})_{i \geq 1}$ over spacetime dyadic points.  It is denoted by $\mathcal{Y}'$.
The second is based on an independent system  of coalescing brownian motions $(B_1^{x_i,t_i})_{i \geq 1}$ and is denoted by $\mathcal{Y}''$. Finally $\mathcal{Y}$ is derived from $(B^{x_i,t_i})$ for $x<0$ and from $(B_1^{x_i,t_i})$ for $x \geq 0$. In the case of coalescence precedence is given to   the  $B_1$ process if the two motions ($B$ and $B_1$) meet in $\R_+ \times \R$ and to  the  $B$ process otherwise.

For an event $A$ we write $A^{\mathcal{Y}}$  to signify that that the event is for the realization $\mathcal{Y}$.
We similarly use the notation  $A^{\mathcal{Y}'}$ and  $A^{\mathcal{Y}''}$.

\bigskip

We now fix $0 < \epsilon <1$ and $0< v_4 < v_2 < A$. By arguing as with Lemma \ref{lem: brown 1} we have that outside probability $K e^{-c N^{\frac{\epsilon}{2}}}$ for some universal    $c =c(v_4,v_2)$ and $K = K(v_4,v_2), \  \exists$ a $B_1^{x_i,t_i}$ with $s_N - N^{-\epsilon} \leq t_i \leq s_N, x_i \in (v_4,v_2)$ so that $B_1^{x_i,t_i}(s) \in (v_4,v_2) \forall s \in (t_i,t_i + 3 N^{-\epsilon})$.  Similarly $\exists B^{x_j,t_j}$ with $x_j \in (-v_2,-v_4), s_N - N^{-\epsilon} \leq t_j \leq s_N$ so that $B^{x_j,t_j}(s) \in (-v_2,-v_4) \forall s \in (t_j, t_j+ 3 N^{-\epsilon})$.

Let $D$ be the event that at least one of these two events fails to occur.  Then  

\begin{eqnarray*}
& &A^{u,\mathcal{Y}}_{a_1,b_1,a_2,b_2}(s_N,t_n) \subseteq  A^{1,\mathcal{Y}}_{a_1,b_1,a_2,b_2}(s_N,t_N) \equiv\\
& &\{ \exists s_N \leq r \leq t_N \hbox{ so that } \exists (\gamma,a), (\gamma',a') \hbox{ with } \sigma= \sigma' = r \hbox{ and } \gamma(\sigma) \in (a_1,b_1), \\ 
& &\gamma'(\sigma') \in (a_2,b_2) \hbox{ on } (\sigma, \sigma+N^{-\epsilon}), a(s) = a'(s) = s-\sigma, \gamma(\sigma+N^{-\epsilon}) \ne \gamma'(\sigma+N^{-\epsilon})\}\\
& \subseteq & D \cup (C_{a_1,b_1}^{\mathcal{Y}'}(s_N,t_N) \cap C_{a_2,b_2}^{\mathcal{Y}''}(s_N,t_N)),
\end{eqnarray*}

where $C_{a_1,b_1}^{\mathcal{Y}'}(s_N,t_N) = \{ \exists \hbox{ path in interval   } (-b_1 -1,-v_4,) 
 \hbox{ of age in interval } [N^{-\epsilon}-1/N, N^{-\epsilon}]\hbox{ at time  } s_N+N^{-\epsilon} \}
\cup  \{  \hbox{ max displacment  of } B^{b_1 - /1/2, s_N} \mbox{ in } (s_N,s_N +N^{-\epsilon} ) \geq 1/2\}$

and similarly for $C_{a_2,b_2}^{\mathcal{Y}''}(s_N,t_N)$.   
By Corollary 2.6 of \cite{MRV}, the probability of $C_{a_1,b_1}^{\mathcal{Y}'}(s_N,t_N)$ is bounded by $ \frac{c}{(N^{- \epsilon} - 1/N)^{1/\alpha}} -\frac{c}{N^{- \epsilon}/ \alpha } $
for constant $c$ depending on $a_2-b_2$ and the constant described in Corollary 2.6 . This is easily bounded by
$c^\prime \frac{1}{N} N^{\epsilon (1+ 1/ \alpha )}$.

But $P(D) \leq k e^{-c N^{\frac{\epsilon}{2}}}$ and so

$P(A^{u,\mathcal{Y}}_{a_1,b_1,a_2,b_2}(s_N,t_n)) \leq K' \frac{N^{\epsilon(2+ 2/ \alpha )}}{N^2} \ < \ K'/N^{3/2}$
for $\epsilon $ chosen sufficiently small.
\end{proof}

\begin{proof}[Proof of Lemma \ref{beegee}]
We begin with some reductions.  First, by brownian scaling it is only necessary to prove the result for $a=0, \ b = 1$.  Secondly, by monotonicity, it is only necessary to show the result for 
$m$ of the form integer part $h^R$ for a fixed $h $ greater than one but sufficiently close.  For notational simplicity we will write $h^R$
for the integer part.  Finally our infinite particle system arises as the limit of finite coalescing  systems.  Accordingly, it is only necessary to treat systems starting with motions at $\Z h^{-M} \cap [a,b]$ provided the bounds obtained are uniform in $M$.

Given $M$, we define the sequence of stopping times for the system $T_r : r \geq 0$ by $T_0 \ = \ 0, $
$$
T_{r+1} \ = \ \inf \{ t \geq T_r: \# \mbox{motions }< h^{M-r-1}\}.
$$
It suffices to prove that for $c $ small (but not depending on $M$ or $r$)
$$
P(T_{r+1} -T_r  > \frac{(1-h^{-2})}{c} h^{-2(M-r)} \mid \mathcal{F}_{T_r} ) \ < \ e^{-ch^{M-r}}
$$
since
$$
P(T_{M-R}> h^{-2R} /c\ \leq \ P(\cup _{j=0}^{M-R-1} (T_{j+1} -T_j > \frac{(1-h^{-2})}{c} h^{-2(M-j)}) \ \leq \ \sum_{k=R+1}^\infty  e^{-ch^{k}}.
$$
But this inequality follows directly from the key argument of \cite{BG} once $h$ is fixed close to $1$ and then $c$ is taken small.
\end{proof}

In order to prove  Lemma \ref{4B} we will need the following lemma.

Recall  probability measure $\nu $ on $\mathbb{R}$ is an $h$-tilting of probability $\mu $ if 
$$
\nu(dx) \ = \ \mu(dx) h(x) / (\int h(y)\mu(dy))
$$
for positive $\mu$-integrable function $h$.

\begin{lemma} \label{tilt}
If $ \nu$ is an $h$ tilt for probability $\mu$ where $h$ is increasing, then 
$$
\nu[x,\infty) \geq (1-c) \mu[x,\infty) \ + \ c
$$
for $c \ = \ \int (h(y)-h(x))_+ \mu (dy) / ( \int (h(y)\mu (dy) )$.
\end{lemma}

{\it Proof of Lemma \ref{4B}} \\

We first note that by a simple FKG argument and scaling, it is enough to show that 
there exists finite $C$ and strictly positive $\epsilon $ so that for
three independent brownian motions $B^1, B^2$ and $B^4$, (the suffix $3$ is reserved for an auxiliary brownian motion $B^3$) with $B^1(0) < B^2(0) < B^4(0) < B^1(0) +3$,  and for all $N$
$$
P( \nexists s \leq 2^N: (B^2(s)-B^1(s))(B^4(s)-B^2(s)) = \  0 ) \leq C2^{-(1+ \epsilon /2)N}.
$$

For such Brownian motions if we condition on the event (of probability of the order $2^{-N/2}$) $\{ \forall 0 \leq s \leq 2^N: B^2(s) > B^1(s)\}$,
then we
expect that motion $B^2$ is stochastically greater than a brownian motion, as it is ``pushed" to the right and so the conditional probability that $\{ \forall 0 \leq s \leq 2^N: B^4(s) > B^2(s)\}$ should be small compared to 
$2^{-N/2}$.

This motivates the  following which is elementary and follows from the analagous result for simple random walks followed by an
application of Donsker's invariance principle:
\begin{lemma}
Given two independent brownian motions $B^1$ and $B^2$with $B^1_0 < B^2_0 $ and conditioned so that
$\forall 0 \leq s \leq T \ B^1_s < B^2_s$ (so that the two brownian motions are no longer  independent and no longer brownian motions), it is possible to define a Brownian motion $B^3$ so that 
\begin{enumerate}
\item
$B^3_0 \ = \ B^2_0 $ and
\item
 for all
$0 \leq s \leq t  \leq T \ B^3_t - B^3_s \ \leq \ B^2_t- B^2_s $.
\end{enumerate}
\end{lemma}

\begin{lemma} \label{bmsandwhich}
For $B^1, B^3, B^2 $ as above  and $K$ and $\epsilon ^\prime $ any fixed values in $(0, \infty )$, there exists $\delta >0$ so that for all positive integer $N$,
outside probability $2^{-3N}$  

$$\sum_{j=1}^{N-1}I_{D^c_j} < N/20
$$ 
for $D_j$ the event that on time interval $\exists \  s  \in \ [2^j,2^{j+1}]$ 
 so that  
\begin{itemize}
\item 
$B^2_s> B^3_s+2 \delta 2^{j/2}$ or
\item
$|B^i_s| \geq \ K 2^{j/2}$ for some $i$ or
\item
$B^2_s - B^1_s \leq \ \epsilon ^\prime 2^{j/2}$
\end{itemize}
\end{lemma} 
\begin{proof}
 It is enough to argue that for any $c > 0$, there exists $\delta >0   $ so that for $j \leq 99/100 N$,
$$
P(A_{j+1} \mid \mathcal{F} _{2^j}) \ \geq \ 1- c
$$
where $ \mathcal{F} _{2^j}$ is the natural filtration of the $B^i$ and $A_{j+1}$ is the event that stopping time $\sigma < 2^{j+1}$ for $\sigma $ the first time after time $2^j$ that \\
$B^2_s -B^3_s \ \geq \ 2 \delta 2^{j/2}$ or\\
$|B^i_s| \geq \ K 2^{j/2}$ or \\
$B^2_s - B^1_s \leq \ \epsilon ^\prime 2^{j/2}$

We fix $\epsilon _1 << \epsilon ^ \prime $ to be specified later and consider $s = 2^{j} +r \epsilon_12^{j}  < \sigma$.  We let $B_1$ evolve on the interval $[s,s+\epsilon_12^{j}]$.  Given this and $B^2_s = z, \ B^2_{s+\epsilon_12^{j}} $ has a tilted disitribution with respect to 
measure $\mu \ = \ g(.-z,\epsilon_12^{j} ) $ for function $h(x)  \ = \ P_{BB}^{(z,s),(x,s+\epsilon_12^{j})}(X$ does not hit $B^1)$ times
$P( \tau \geq 2^N \mid B^a([0,s+\epsilon_12^{j})], B^2(s+\epsilon_12^{j}) = x)$.  
Here $P_{BB}$ denotes the probability for space time Brownian bridges.  We note that both these functions are increasing in $x$.
Furthermore unless $\sigma$ hs aleady occurred  the function $$P( \tau \geq 2^N \mid B^1([0,s+\epsilon_12^{j}]), B^2(s+\epsilon_12^{j}) = x)$$ is at least linear with slope at least $\frac{1}{2 2^{N/2}(2K+ \epsilon ^\prime )}$ for $N$ large within the relevent interval.  From this and Lemma \ref{tilt} we obtain that 

$P(B^2_{s+\epsilon_12^{j}} >B^3_{s+\epsilon_12^{j}} +2 \delta 2^{j/2}   )$ with probability $C(K, \epsilon^\prime )$ for $\epsilon$ small.  
So 
$$
P(\sigma > 2^j + r \epsilon_1 2^j) \ < \ (1-C)^r.
$$
By taking $\epsilon _1$ sufficiently small we obtain the desired result.

\end{proof}
By property (2) of brownian motion $B^3$ we have
\begin{corollary} \label{corbmsandwhich}
For $\delta $ as above given $K$ and $\epsilon ^\prime $,
Let $D_j ^\prime$ be the event that 
\begin{itemize}
\item 
$ \forall s \in (2^j, 2^{j+1} )B^2_s> B^3_s+ \delta 2^{j/2}$ or
\item
$|B^i_s| \geq \ K 2^{(j-1)/2}$ for some $i$ and some $s \in (2^{j-1} ,2^j)$ or
\item
$B^2_s - B^1_s \leq \ \epsilon ^\prime 2^{(j-1)/2}$  for some $s \in (2^{j-1} ,2^j) $
\end{itemize}

Then outside probability $2^{-3N}$  
$\sum_{j=1}^{N-1}I_{(D^\prime_j)^c} < N/20$
\end{corollary}
The proof of Lemma \ref{4B} (which follows the path of \cite{Kes}) consists largely of making a list of ``bad" events $A^i_j, \ i=1,2 \cdots 5$ and $1 \leq j \leq N$ which concern the motions 
$B^1, B^2 $ and $B^3$ on time interval $[2^j, 2 ^{j+1}]$ and such that 
$$
\sum_{j=1}^N I_{\cup_{i=1} ^5 A^i_j} \ \leq \ N/4
$$
outside probabilty $C2 ^{-3N}$.  We then argue that conditioned on a brownian motion $B^4$ starting to the right of $B^3$
avoiding it until time $2^N$, the probability that it avoids $B^2$  is still  exponentially  small in $N$.

We first pick (see \cite{Kes}) $K$ so large for event
$$
A_j^1 \ = \ \{  \exists t \in [2^j, 2^{j+1}], i=1,2 \mbox{ or } 3: |B^i_t| \geq K2^{j/2}\}
$$
$\{ \sum_{j=1}^N I_{ A^1_j} \ \geq \ N/100\}$ has probabilty  less than $2^{-3N}$ for all $N$.

Let $A^2_j$ be the event that for some $i = 1,2$ or $3$ $T^i_j \notin ((1+\epsilon_0 )2^j, (1-\epsilon_0)2^{j+1})$,
where $T^i_j$ is the (a.s. unique) time in $(2^j, 2^{j+1})$ at which brownian motion $B^i $ achieves its maximum on this time interval.  It is
clear that we can choose $\epsilon_0$ so small that $\{ \sum_{j=1}^N I_{ A^2_j} \ \geq \ N/100\}$ has probabilty  less than $2^{-3N}$ for all $N$.

 $A^3_j$ is the event that for some $t \in  [2^j, 2 ^{j+1}], B^2_t - B^1_t \ \leq \ \epsilon^\prime  2^j$.  Again (this time using the conditioning event defining $B^1$ and $B^2$) we have the desired conclusion for $\epsilon^\prime $ fixed sufficiently small.

$A^4_j$ is the complement of the event that on interval $[2^j, 2 ^{j+1}], \ B^2_t - B^3_t \geq \delta2^j \ \forall t$. Here $\delta $ is chosen as in Lemma \ref{bmsandwhich} and Corollary \ref{corbmsandwhich}, given $K$ and $\epsilon ^\prime $ specified in defining the events $A^1_j $ and $A^3_j$.
We have  by  Corollary \ref{corbmsandwhich} and our choices of $K$ and $\epsilon ^\prime $ the desired inequality   $P (\sum_{j=1}^N I_{ A^4_j} \ \geq \ N/10) | \leq \ 2^{-3N}$  for all $N$.

Finally $A^5_j$ is the  complement of the  event that for $\ \epsilon _3>0$, we have 
$|B^i_s - B^i_{2^j}| \leq \ \epsilon^\prime  2^{j/2}/ 10 \forall s \in [2^j, 2^j (1+ \epsilon _3) )$
and $i=1,2$ and 
$|B^2_s - B^2_{2^{j+1}}| \leq \ \epsilon ^\prime  2^{j/2} / 10\forall s \in [2^{j+1}(1- \epsilon _3 ), 2^{j+1}) $
.Again if we choose $\epsilon _ 3 $ sufficiently small then we have for all $N$,  outside probability $2^{-3N}$, that

$$
\sum_{j=1}^N I_{ A^5_j} \ \leq \ N/10) 
$$

Thus we have that outside of probability $C 2^{-5N/2}, \ \sum_{j=1}^N I_{\cup_{i=1} ^5 A^i_j} \ \leq \ N/4$.

Given the realizations of the  $B^i: \ i \leq 3$ we choose $B^4$ conditioned to not meet $B^3$ (which is necessary if $B^4$is to avoid $B^2$).   
This conditioning event has probability bounded by $C2^{-N/2}$ as is well known.  Given this we condition on the values of $B^4(2^j)$ so that on interval $(2^j,2^{j+1})$, $B^4$ is a Brownian bridge conditioned to avoid $B^3$,  conditionally independent of $B^4 $ outside this interval given .$B^4(2^j)$ and $B^4(2^{j+1})$. 

By our choice of $K$ and $ \epsilon ^\prime $, we have that outside of probability $C2^{-5N/2}$ for some $C$,
the number of $j \leq N $ satisfying 
\begin{itemize}
\item
$|B^4(2^j)| $ or $|B^4(2^{j+1})|  \ \geq \ K 2^{j/2}$ or
\item
$ B^4(2^j) -  B^3(2^j)\ < \ \epsilon^\prime  2^{j/2}$ or 
\item
$ B^4(2^{j+1}) -  B^3(2^{j+1})\ < \ \epsilon ^\prime  2^{j/2}$ 
\end{itemize}
is less than $N/4$.
A $j$ for which none of the above is true and for which none of the $A^j$ occur is said to be {\it good}.
We have that outside of probability $C2^{-5N/2}$ for some universal $C$, the number of good j exceeds $N/2$.

But it is easily seen that (by scaling), there exists $\delta ^\prime \ = \  \delta ^\prime (K, \epsilon ^\prime, \delta, \epsilon_2, \epsilon _3) >0$ so that for any good $j$ the probability that $B^4$ hits $B^2$ (conditioned upon not hitting $B^3$ ) is at least $\delta ^\prime$ independently of the event for other values of $j$.  This suffices.

\end{document}